\documentclass{article}
\usepackage{times}
\usepackage[margin=1.3in]{geometry}
\usepackage{helvet}  
\usepackage{courier}  
\usepackage{url}  
\usepackage{graphicx}  
\usepackage[utf8]{inputenc} 
\usepackage{booktabs}       
\usepackage{amsfonts}       
\usepackage{nicefrac}       
\usepackage{microtype}      
\usepackage{amsmath}
\usepackage{amssymb}
\usepackage{amsthm}
\usepackage{algorithm}
\usepackage{algorithmic}
\usepackage{subfigure}
\usepackage[usenames, dvipsnames]{color}
\usepackage{flexisym}

\DeclareMathOperator*{\argmin}{arg\ min}   
\DeclareMathOperator*{\argmax}{arg\ max}   
\DeclareMathOperator*{\prox}{prox}
\DeclareMathOperator*{\Prog}{Prog}
\DeclareMathOperator*{\dom}{dom}
\DeclareMathOperator*{\with-probability}{with\ probability}
\DeclareMathOperator*{\otherwise}{otherwise}

\newcommand*{\QEDB}{\hfill\ensuremath{\square}}

\newtheorem{assumption}{\textbf{Assumption}}
\newtheorem{theorem}{\textbf{Theorem}}
\newtheorem{lemma}{\textbf{Lemma}}
\newtheorem{corollary}{\textbf{Corollary}}

\newcommand{\tabincell}[2]{\begin{tabular}{@{}#1@{}}#2\end{tabular}}

\title{Exploring Fast and Communication-Efficient Algorithms in Large-scale Distributed Networks}
\author{Yue Yu$^{1}$, Jiaxiang Wu$^{2}$, Junzhou Huang$^{2,3}$\\
$^1$Institute for Interdisciplinary Information Sciences, Tsinghua University\\
$^2$Tencent AI Lab\\
$^3$University of Texas at Arlington\\
yu-y14@mails.tsinghua.edu.cn,
jonathanwu@tencent.com,
jzhuang@uta.edu
}

\date{}

\begin{document}

\maketitle

\begin{abstract}
The communication overhead  has become a significant bottleneck in data-parallel network with the increasing of model size and data samples.
In this work, we propose a new algorithm LPC-SVRG with quantized gradients
and its acceleration ALPC-SVRG to effectively reduce the communication complexity while maintaining the same convergence as the unquantized algorithms.
Specifically, we formulate the heuristic gradient clipping technique within the quantization scheme and
show that  unbiased quantization methods  in related works \cite{alistarh2017qsgd,wen2017terngrad,zhang2017zipml} are special cases of ours.
We introduce \emph{double sampling} in the accelerated algorithm ALPC-SVRG to fully combine the gradients of full-precision and low-precision, and then achieve acceleration with fewer communication overhead.
Our analysis focuses on the nonsmooth composite problem, which makes our algorithms more general.
The experiments on linear models and deep neural networks validate the effectiveness of our algorithms.
\end{abstract}

\section{INTRODUCTION}
Recent years has witnessed data explosion and increasing model complexity in machine learning.
It becomes difficult to handle all the massive data  and large-scale models within one machine.
Therefore, large-scale distributed optimization receives growing attention \cite{dean2012large,abadi2016tensorflow,de2015taming,xing2015petuum}.
By exploiting multiple workers, it can remarkably reduce computation time.
%
Distributed data-parallel network is a commonly used large-scale framework which contains $N$ workers
and each worker keeps a copy of model parameters.
At each iteration, all workers compute their local gradients and communicate gradients with peers to obtain global gradients, and update model parameters.
As the number of workers increases, the computation time (for a mini-batch of the same size) can be dramatically reduced, however, the communication cost rises.
It has been observed in many distributed learning systems that the communication cost has become the performance bottleneck \cite{chilimbi2014project,seide20141,strom2015scalable}.
\\
\\
\noindent
To improve the communication efficiency in data-parallel network, generally, there are two orthogonal methods, i.e., gradient quantization and gradient sparsification.
Researches on quantization focus on employing low-precision and fewer bits representation of gradients rather than full-precision with $32$ bits\footnote{In this paper, we assume that a floating-point number is stored using 32 bits even though it can be 64 bits in many modern systems. }~\cite{seide20141,wu2018error,zhang2017zipml,wen2017terngrad,alistarh2017qsgd}.
And for works on sparsification, they design dropping out mechanisms for gradients to reduce the communication complexity \cite{wangni2017gradient,aji2017sparse}.
\\
\\
\noindent
%
For most gradient quantization based methods, the unbiased stochastic quantization is adopted to compress gradients into their low-precision counterparts.
However, it has been observed in many applications that such unbiased quantization brings larger precision loss or quantization error \cite{wen2017terngrad},
and therefore leads to lower accuracy.
%
%
On the other hand, several works \cite{wen2017terngrad,kanai2017preventing,pascanu2013difficulty} reported the heuristic \emph{gradient clipping} technique
can effectively improve convergence.
However, this breaks the unbiasedness property and no theoretical analysis has been given so far.
\\
\\
\noindent
In our work,
we propose the LPC-SVRG algorithm to embed gradient quantization and clipping techniques into SVRG~\cite{johnson2013accelerating}. Furthermore, we present its acceleration variant, the ALPC-SVRG algorithm.
These two algorithms
adopt the variance reduction idea to reduce the gradient variance as the algorithm converges, so as to achieve a faster convergence rate.
%
To reduce the communication overhead, we introduce a new quantization scheme which integrates gradient clipping, and provide its theoretical analysis.
%
%
%
%
%
Furthermore, we propose 
%
%
%
\emph{double sampling} to fully take advantage of both low-precision and full-precision representations,
and achieve lower communication complexity and fast convergence rate at the same time.
Detailed contributions are summarized as follows.
\\
\\
\noindent
\textbf{Contribution.}
We consider the following finite-sum composite minimization problem:
\begin{equation}
  \label{eq:pro-formu}
  \min_{x \in \mathbb{R}^d} P(x) = f(x) + h(x) = \frac{1}{n}\sum\limits_{i=1}^n f_i(x) + h(x),
\end{equation}
\noindent
where $f(x)$ is an average of $n$ smooth functions $f_i(x)$,
and $h(x)$ is a convex but can be nonsmooth function.
$d$ is the model dimension.
Such formulation generalizes many applications in machine learning and statistics, such as logistic regression and deep learning models.
\\
\\
\noindent
We try to solve (\ref{eq:pro-formu}) in a data-parallel network with  $N$ workers,
and analyze the convergence behavior when low-precision quantized gradients are adopted.
Our analysis covers both convex and nonconvex objectives, specifically:
\\
\\
\noindent
$(1)$ We integrate gradient clipping into quantization to reduce the overhead of communication,
%
and show that the unbiased quantization methods adopted in \cite{alistarh2017qsgd,wen2017terngrad,zhang2017zipml} are special cases of our scheme.
We also mathematically analyze the quantization error given gradient clipping;
\\
\\
\noindent
$(2)$ Based on such quantization scheme, we propose LPC-SVRG to solve nonconvex and nonsmooth objective (\ref{eq:pro-formu}).
We prove the same convergence rate can be achieved even the numbers are represented by much fewer $O(\log{\sqrt{d}})$ bits (compared to $32$);
\\
\\
\noindent
$(3)$ We propose an accelerated algorithm ALPC-SVRG based on Katyusha momentum \cite{allen2017katyusha}.
 With \emph{double sampling}, we are able to combine updates both from quantized gradients and full-precision gradients to
  achieve fast convergence without increasing the communication overhead.
  \\
  \\
  \noindent
$(4)$ We conduct extensive experiments on both linear regression and deep learning models to validate the effectiveness of our methods.

\section{RELATED WORK}
Many literatures focus on designing fast and efficient algorithms for large-scale distributed systems \cite{abadi2016tensorflow,bekkerman2011scaling,de2015taming,iandola2016firecaffe}.
Asynchronous algorithms with stale gradients are also extensively studied \cite{recht2011hogwild,lian2015asynchronous,dean2012large,li2014communication}, which are orthogonal to our work.
Researches on reducing communication complexity in data-parallel network can be generally divided into two categories: gradient sparsification and gradient quantization.
\\
\\
\noindent
\textbf{Gradient sparsification.}
\cite{aji2017sparse} introduced a dropping out mechanism where only gradients exceed a threshold being transmitted.
\cite{wangni2017gradient} formulated the sparsification problem into a convex optimization to minimize the gradient variance.
Several heuristic methods such as momentum correction, local gradient clipping and momentum factor masking were adopted in \cite{lin2017deep} to compensate the error induced by gradient sparsification.
\cite{wang2018atomo} proposed a method for atomic sparsification of gradients while minimizing variance.
\\
\\
\begin{figure}
\begin{minipage}{0.5\linewidth}
  \centering
  \label{fig:comm-network}
  \includegraphics[width=2in]{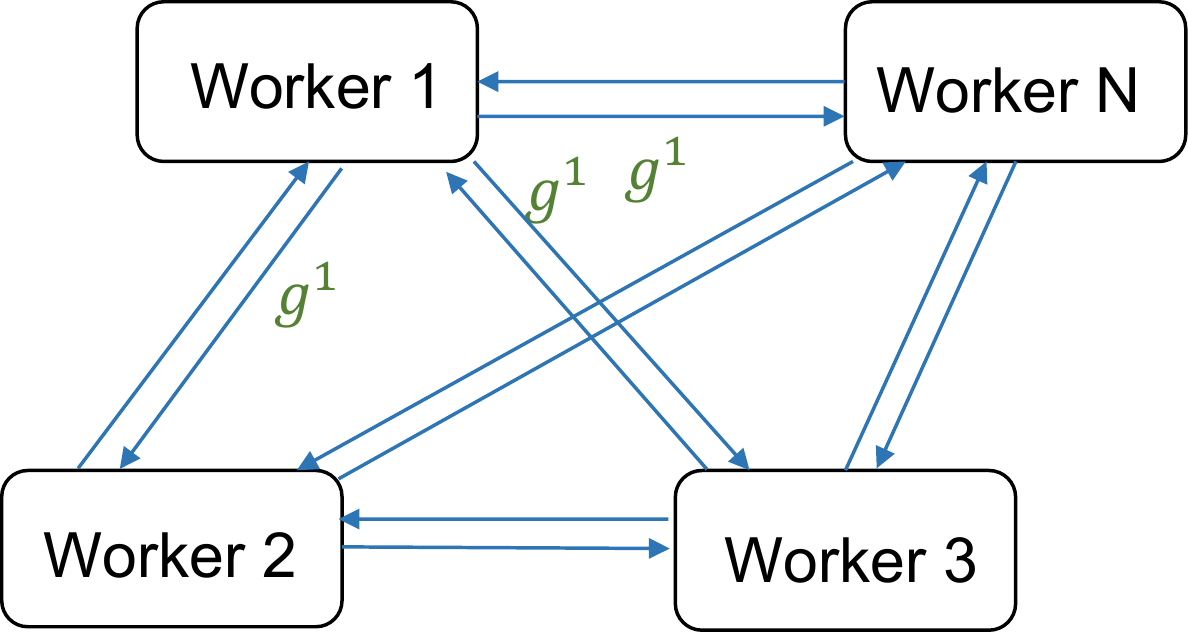}
\end{minipage}%
\begin{minipage}{0.5\linewidth}
  \centering
  \includegraphics[width=2in]{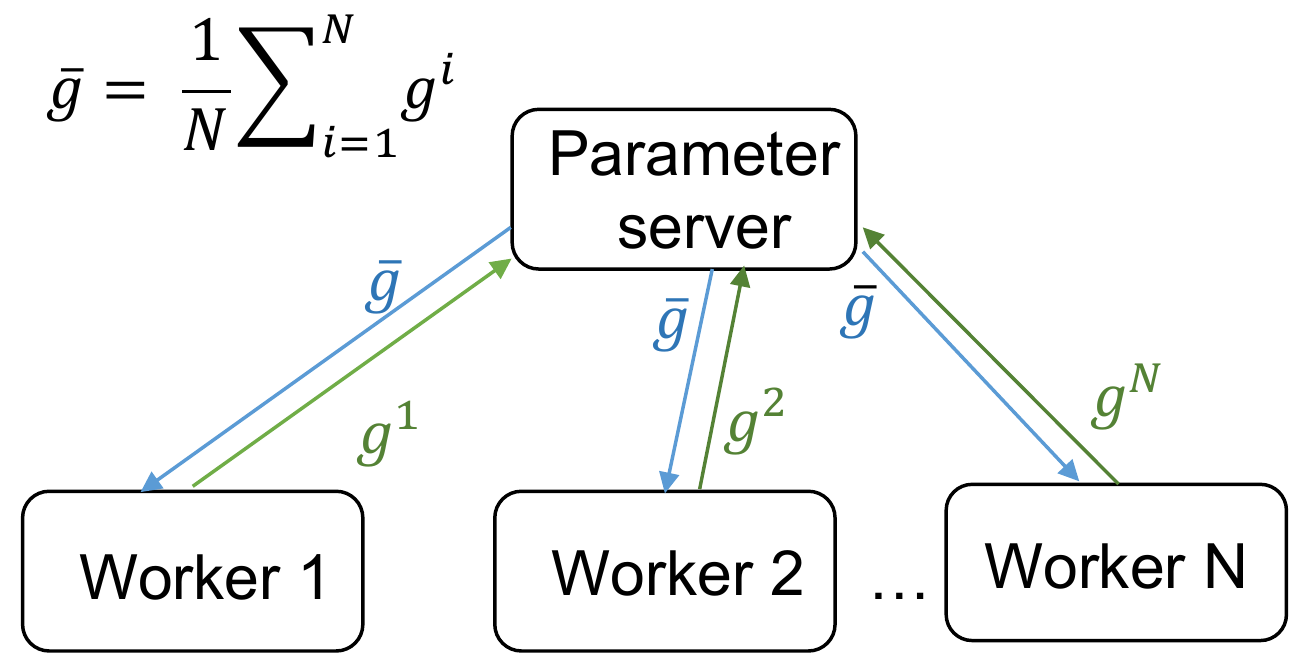}
\end{minipage}
\caption{An illustration of two commonly used communication frameworks: broadcast (left) and parameter-server (right).}
\end{figure}
\noindent
\textbf{Gradient quantization.}
\cite{seide20141} quantized gradients to \{-1,1\} using zero-thresholding and used an error-feedback scheme to compensate the quantization error.
Similar technique was also adopted in \cite{wu2018error}, where original unbiased quantization method was used to quantize the error-compensated gradients.
\cite{wen2017terngrad} adopted a $3$-level representation, i.e. \{-1, 0, 1\}, for each element of gradient.
They empirically showed the accuracy improvement when gradient clipping was applied, while no theoretical analysis was provided.
\cite{alistarh2017qsgd} used an unbiased stochastic quantization method to quantize gradients into s-level.
%
%
They provided the convergence property of unbiased low-precision SVRG \cite{johnson2013accelerating}
for strongly convex and smooth problems.
\cite{de2018high} proposed a low-precision variance reduction method for training in a single machine, with quantized model parameters.
An end-to-end low-precision mechanism was also setup in \cite{zhang2017zipml,zhou2016dorefa}, where data samples, models and gradients were all quantized.
\cite{zhang2017zipml} analyzed the convergence property of an unbiased stochastic quantization for gradients based on a convex problem,
and designed an optimal quantization method for data samples.
\\
\\
\noindent
Our work distinguishes itself from the above results in:
(1) considering general nonsmooth composite problems and providing convergence analysis in both convex and nonconvex settings;
(2) providing theoretical understandings of a new quantization scheme with gradient clipping;
(3) proposing \emph{double sampling} in accelerated algorithms to maintaining both fast convergence and lower communication overhead.
%
\section{PRELIMINARY}
\subsection{Data-Parallel Network}
\label{sec:data-parallel-ar}
We consider a data-parallel distributed network with $N$ workers.
Each worker maintains a local copy of estimation model and can get access to the whole datasets stored in e.g., HDFS.
At each iteration, every worker computes local stochastic gradient on a randomly sampled mini-match.
Then they synchronously communicate gradients with peers to compute global gradients and update models.
%
%
There are two commonly used communication frameworks, as shown in Figure $1$,
to realize the above data-parallel network,
i.e.,
(1) \textbf{broadcast:}
each worker broadcasts local gradients to all peers.
When one worker gathers all gradients from other workers, it averages them to obtain global gradients;
(2) \textbf{parameter-server}\footnote{We still use the convention ``parameter-server'' even though only gradients are transmitted.}\textbf{:}
%
%
%
%
local gradients computed by all workers are synchronously gathered in the parameter server.
\\
\\
\noindent
%
Selecting a proper communication framework is determined by real-world requirements.
In Section \ref{sec:Communication Schemes and Complexities}, we provide the communication overhead comparisons 
when different frameworks are plugged into our algorithms.
%
\subsection{Low-Precision Representation and Quantization}
%
In this paper, we adopt a low-precision representation, denoted as a tuple $(\delta,b)$, of transmitted numbers.
%
%
The procedure of transforming numbers from full-precision to low-precision is called quantization.
\\
\\
\noindent
%
For tuple $(\delta,b)$, $b \in \mathbf{N}$ represents the amount of bits (used for storing numbers) and $\delta \in \mathbb{R}$ is the scale factor,
its representation domain or quantization codebook is
\begin{equation*}
  \dom(\delta,b) = \{-2^{b-1}\cdot\delta,..., -\delta, 0, \delta,..., (2^{b-1}-1)\cdot \delta\}.
\end{equation*}
Denote $Q_{(\delta,b)}(x)$ as the quantization function for a number $x \in \mathbb{R}$ given $(\delta,b)$.
It outputs a point within $\dom(\delta,b)$ in the following rules:
\\
\\
\noindent
%
%
(1) if $x$ is in the convex hull of $\dom(\delta, b)$, without loss of generality, $x \in [z, z+\delta]$, where $z \in \dom(\delta,b)$.
It will be stochastically rounded up or down in the following sense:
\begin{equation*}
  Q_{(\delta,b)}(x) =
  \begin{cases}
    z, \quad \quad \ \ \with-probability \ \frac{z+\delta-x}{\delta},\\
    z+\delta, \quad   \otherwise.
  \end{cases}
\end{equation*}
It can be verified that the above quantization is unbiased, i.e., $\mathbf{E}[Q_{(\delta,b)}(x)] = x$.
%
%
%
Moreover, the quantized variance or precision loss can be bounded as:
\begin{lemma}
\label{lemma:unbia-variance}
  If $x$ is in the convex hull of $\dom(\delta,b)$, then we have
  \begin{equation*}
  \mathbf{E}[\big( Q_{(\delta,b)}(x) - x \big)^2] \leq \delta^2/4.
\end{equation*}
\end{lemma}
\noindent
The intuition behind Lemma~\ref{lemma:unbia-variance} is that a smaller $\delta$ determines more dense quantization points, and as a result incurs less precision loss (or quantization error).
\\
\\
\noindent
(2) on the other hand, if $x$ is not in the convex hull of $\dom(\delta,b)$, it will be projected to the closest point, in other words, the smallest or largest value in $\dom(\delta,b)$.
\\
\\
\noindent
In the following sections, we use function $Q_{(\delta,b)}(\cdot)$ to quantize gradient, which means each coordinate is quantized using the same scale factor independently.
Also note that low-precision numbers with the same scale factor can be easily added with several bits overflow.

\begin{algorithm}[!tb]
   \caption{Low-Precision SVRG with Gradient Clipping: \textbf{LPC-SVRG} (for each worker $i$)}
\begin{algorithmic}[1]   \label{Algm:LPC-SVRG}
   \STATE {\bfseries Input:} $S$, $m$, $\lambda$, $B$,  $\eta$, $\tilde{x}^0 = x^0$;
   \FOR{$s=0,1,...,S-1$}
       \STATE $x_0^{s+1} = \tilde{x}^s$;
       \STATE \textbf{Communication step $1$:} cooperates with other workers to compute $\nabla f(\tilde{x}^s)$;
       \FOR{$t=0$ {\bfseries to} $m-1$}
           \STATE  uniformly and independently samples (with replacement) a mini-batch $i_B$ with size B to calculate\\
            $u_{t}^{i} = \frac{1}{B} \sum\limits_{a \in i_B} \big[ \nabla f_a(x_t^{s+1}) - \nabla f_a(\tilde{x}^s) \big]$;
           \STATE  \textbf{Quantization step:} $\tilde{u}_{t}^i = Q_{(\delta_{t}^i, b)} (u_{t}^i)$;
           \STATE  
           \textbf{Communication step $2$:} cooperates with other workers to compute $\tilde{u}_t$ using equation (\ref{eq:tilde_u_t}) or (\ref{eq:re-quan});
           \STATE $v_{t}^{s+1} = \tilde{u}_t + \nabla f(\tilde{x}^s)$;
           \STATE $x_{t+1}^{s+1} = \prox_{\eta h}(x_t^{s+1} - \eta v_t^{s+1})$;
       \ENDFOR
      \STATE $\tilde{x}^{s+1} = x_m^{s+1}$;
  \ENDFOR
   \STATE {\bfseries Output:} Uniformly choosing from $\{\{x_t^{s+1}\}_{t=0}^{m-1}\}_{s=0}^{S-1}$.
\end{algorithmic}
\end{algorithm}
\subsection{Notation}
Throughout the paper, we denote $x^*$ as the optimal solution of (\ref{eq:pro-formu}).
$d$ is the dimension of $x$.
%
%
$||\cdot||_{\infty}$ represents the max norm of a vector.
$||\cdot||$ denotes $L_2$ norm and the base of the logarithmic function is $2$ if without special annotation.
\\
\\
\noindent
\textbf{Proximal operator.}
To handle the nonsmooth $h(x)$ in $(\ref{eq:pro-formu})$, we apply the proximal operator which is formed as
$\prox_{\eta h}(x) = \argmin_{y}(h(y) + \frac{1}{2\eta}||y-x||^2)$, where $\eta > 0$.
Proximal operator can be seen as a generalization of projection. If $h(x)$ is an indicator function of a closed convex set, the proximal operator becomes a projection.
\\
\\
\noindent
\textbf{Convergence metric.}
%
%
For convex problem, we simply use the objective gap, i.e., $P(x) - P(x^*)$ as the measurement of convergence,
while for nonconvex nonsmooth problem, we adopt a commonly used metric \emph{gradient mapping} \cite{nesterov2013introductory}:
 $G_{\eta}(x) \triangleq \frac{1}{\eta}[ x- \prox_{\eta h}(x-\eta\nabla f(x))]$.
A point $x$ is defined as an $\epsilon$-accurate solution if $\mathbf{E}||G_{\eta}(x)||^2 \leq \epsilon$ \cite{Reddi2016Fast, Sra2012Scalable}.

\subsection{Assumption}
We state the assumptions made in this paper, which
are mild and are often assumed in the literature, e.g.,  \cite{johnson2013accelerating,Reddi2016Fast}.
\begin{assumption}
  \label{assum:unbia}
The stochastic gradient is unbiased, i.e.,  for a random sample $i \in \{1,...,n\}$,
  $ \mathbf{E}_{i}\big[ \nabla f_{i}(x) \big] = \nabla f(x)$.
Moreover, the random variables sampled in different iterations are independent.
\end{assumption}
\begin{assumption}
  \label{assum:lip}
We require each function $f_i(x)$ in (\ref{eq:pro-formu}) is $L$-(Lipschitz) smooth, i.e.,
$||\nabla f_i(x) - \nabla f_i(y)|| \leq L||x-y||$.
\end{assumption}

\section{LOW-PRECISION SVRG WITH GRADIENT CLIPPING}
Now we are ready to introduce our new  LPC-SVRG algorithm with low-precision quantized gradients.
 As shown in Algorithm~\ref{Algm:LPC-SVRG},
%
%
%
LPC-SVRG divides the optimization process into epochs, similar to
full-precision
SVRG~\cite{johnson2013accelerating,xiao2014proximal,Reddi2016Fast}, and each epoch contains $m$ inner iterations.
At the beginning of each epoch, we keep track of a full-batch gradient $\nabla f(\tilde{x})$, where $\tilde{x}$ is a reference point and is updated at the end of the last epoch.
The model $x$ is updated in the inner iteration, with a new stochastic gradient formed in Steps $6$-$9$.
Here $\tilde{u}_t$ is an averaged quantized gradient and will be specified later based on different communication schemes.
We consider a mini-batch setting, with each node $i$ uniformly and independently sampling a mini-batch $i_B$ to compute stochastic gradient.
%
If no gradient quantization is applied,
 i.e., $\tilde{u}_t = \frac{1}{N}\sum\limits_{i=1}^N u_t^i$ in Algorithm~\ref{Algm:LPC-SVRG}, \cite{Reddi2016Fast} has shown:
\begin{theorem}
  \label{thm:full-precision-svrg}
(\cite{Reddi2016Fast}, Theorem $5$). Suppose Assumptions~\ref{assum:unbia}, \ref{assum:lip} hold, $h(x)$ is convex and $NB \leq n$. Let $T$ be a multiple of $m$ and $\eta = \rho/L$ where $\rho < 1/2$ and satisfies
  $4\rho^2m^2/(NB) + \rho \leq 1$.
Then for the output $x_{out}$ of Algorithm $1$ (with $\tilde{u}_t = \frac{1}{N}\sum\limits_{i=1}^N u_t^i$), we have:
\begin{equation}
  \mathbf{E}[||G_\eta(x_{out})||^2] \leq \frac{2L(P(x^0) - P(x^*))}{\rho(1-2\rho)T}.
\end{equation}
\end{theorem}

\subsection{Low-Precision with Gradient Clipping}
\label{sec:Low-Precision with Gradient Clipping}
As the scale of data-parallel network and model parameters increases, the communication overhead greatly enlarges and even becomes the system performance bottleneck \cite{zhang2017zipml,seide20141,strom2015scalable}.
In Algorithm~\ref{Algm:LPC-SVRG},
we propose a gradient quantization
 method
to reduce communication complexity.
%
%
Specifically, in Step $7$, we independently quantize each coordinate of $u_t^i$ using tuple $(\delta_t^i,b)$ before it is transmitted.
The values of $\delta_t^i$ and $b$ determine the precision loss and the amount of communication bits.
In our work, the scale factor $\delta_t^i$ is set to be
\begin{equation}
  \label{eq:delta_t_i}
  \delta_t^i = \frac{\lambda ||u_t^i||_{\infty}}{2^{b-1}-1},
\end{equation}
where $\lambda \in (0,1]$ plays the role of clipping parameter.
\\
\\
\noindent
\textbf{Unbiased quantization.}
If $\lambda = 1$,
it can be verified that each coordinate of $u_t^i$ is in the convex hull of $\dom(\delta_t^i,b)$.
Thus, $\mathbf{E}[\tilde{u}_t^i] = u_t^i$.
Such unbiased quantization method is equivalent to the quantization function adopted in
 \cite{alistarh2017qsgd,wen2017terngrad,zhang2017zipml,wu2018error},
%
where the number of quantization levels adopted equals to the number of positive points in $\dom(\delta,b)$.
%
\\
\\
\noindent
\textbf{Biased quantization with gradient clipping.}
Lemma~\ref{lemma:unbia-variance} shows that a larger $\delta_t^i$ leads to a bigger quantization error.
Thus, for each coordinate in the convex hull of $\dom(\delta_t^i,b)$,
%
the quantization error will be reduced by a factor of $\lambda^2$ if $\lambda < 1$.
On the other hand, when $\lambda < 1$,
%
%
there exist coordinates exceeding $\dom(\delta_t^i,b)$, which will be projected to the nearest representable values.
%
That's we call gradient clipping.
%
%
%
%
\begin{figure}
  \label{fig:gaussian}
\begin{minipage}{0.5\linewidth}
  \centering
  \includegraphics[width=2in]{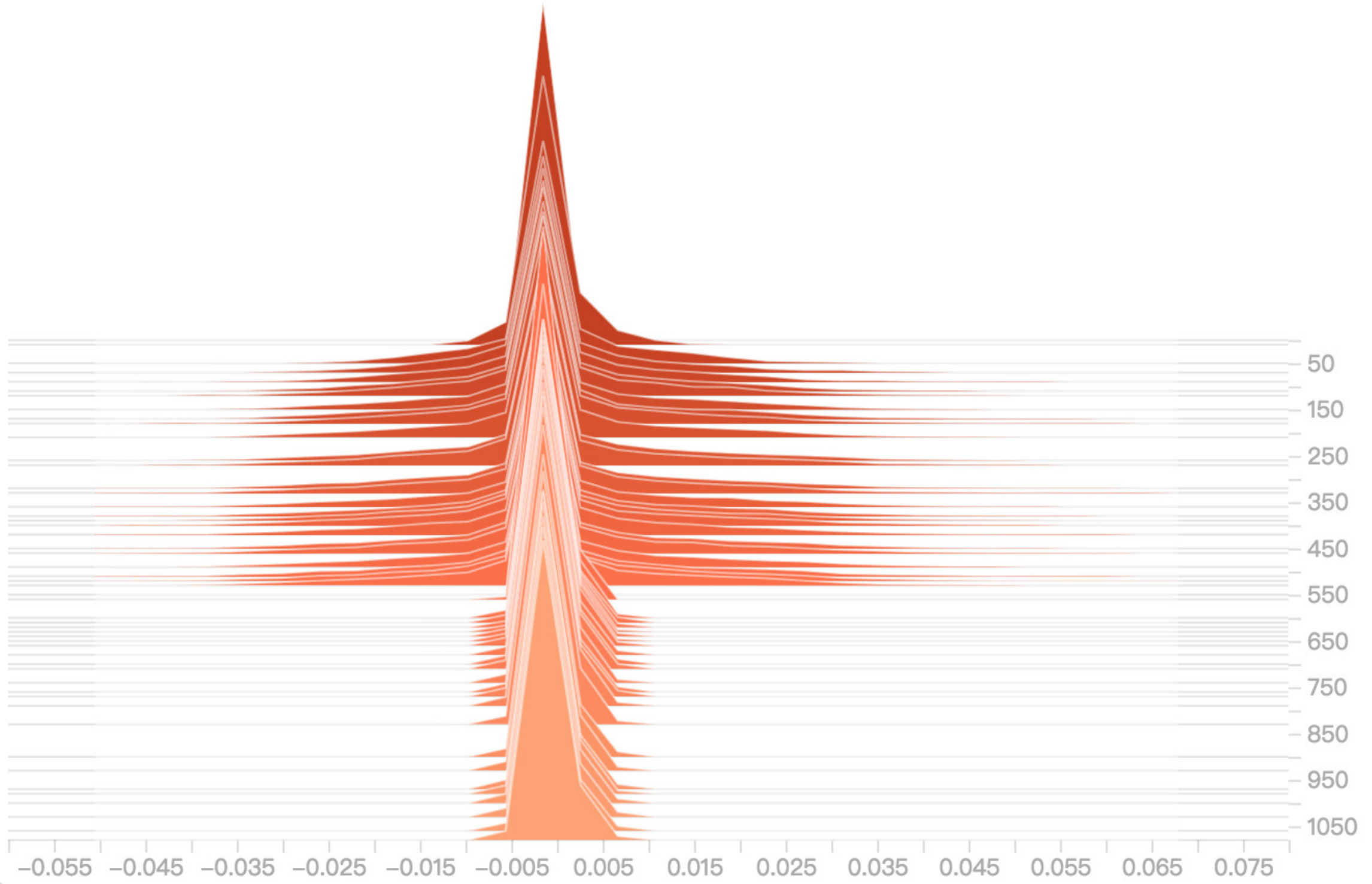}
\end{minipage}%
\begin{minipage}{0.5\linewidth}
  \centering
  \includegraphics[width=2in]{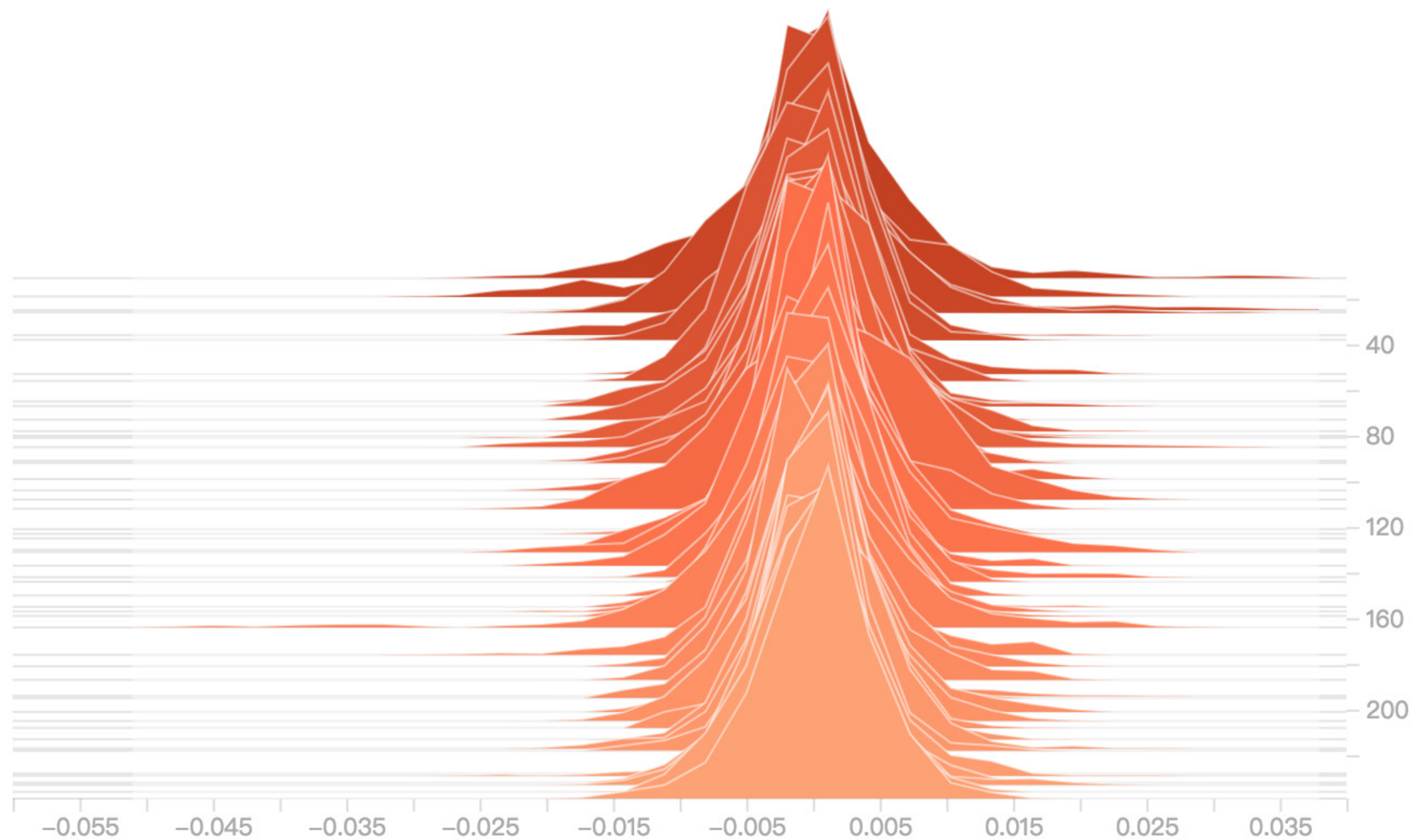}
\end{minipage}
\caption{Histograms of gradients $u_t^i$ (see Algorithm~\ref{Algm:LPC-SVRG}),
left: logistic regression for  MNIST; right: the second convolutional layer of ResNet-20 for CIFAR-10.
The vertical axis represents the index of training iterations.}
\end{figure}
\\
\\
\noindent
To demonstrate the intuition of gradient clipping, in Figure 2, we visualize
the histograms of $u_t^i$ calculated in ResNet-20 model \cite{he2016identity} on dataset CIFAR-10  \cite{krizhevsky2009learning},
and logistic regression on  MNIST \cite{chang2011libsvm}.
We can see that the distribution of $u_t^i$ is very close to a Gaussian distribution and
very few coordinates exceed region $[-c\sigma, c\sigma]$, where $\sigma$ is the approximated standard deviation.
Therefore, setting $\lambda=1$ incurs great quantization error since $||u_t^i||_{\infty}$ is usually far more large.
\\
\\
\noindent
Such Gaussian distribution (of gradients in SGD) has also been recorded in \cite{wen2017terngrad}, and they conducted cross validation to choose the clipping parameter.
However, only convergence proof under unbiased quantization with no gradient clipping was provided.
%
%
In this work, we formally integrate gradient clipping technique into quantization through a parameter $\lambda$ and
 %
provide  the theoretical analysis for biased quantization in Section~\ref{sec:LPC-SVRG-ana}.

 \begin{table*}[h]
 \caption{Comparisons of different schemes for \textbf{Communication step $2$} in Algorithm~\ref{Algm:LPC-SVRG}.
 The values of the last column equal to the ratio of bits used by full-precision transmission to that of low-precision,
 e.g., for \textbf{Broadcast} scheme, it equals to $32dN(N-1)/(32+bd)N(N-1)$.}
 \label{table:com-scheme}
 \begin{center}
 \begin{small}
 \begin{tabular}{ccccc}
 \toprule
 Communication scheme & With re-quantizing & Bits overflow & \# of transmitted bits & Communication reduced factor \\
 \midrule
 Broadcast   & No & No & $(32+bd)N(N-1)$ & $\frac{32}{32/d + b}$\\
 Parameter-server   & No & Yes & $ N(64+2bd+d\lceil\log{N}\rceil)$ & $\frac{32}{32/d + b + (\lceil\log{N}\rceil)/2}$\\
 \tabincell{c}{Parameter-server}  & Yes & No & $N(64+2bd)$ & $\frac{32}{32/d + b}$ \\
 \bottomrule
 \end{tabular}
 \end{small}
 \end{center}
 \end{table*}
\subsection{Communication Schemes and Complexities}
\label{sec:Communication Schemes and Complexities}
%
In this section, we present the communication schemes for the above quantized low-precision gradients.
As shown in Algorithm~\ref{Algm:LPC-SVRG}, there are two communication steps.
%
%
The first one only involves full-precision operations, and its communication overhead can be neglected since it only happens once per epoch.
%
Below, we provide more details for \textbf{Communication step $2$}.
First note that when the $i$-th worker transmits $\tilde{u}_t^i$,  it only needs to send b bits representations of all coordinates and an extra float $\delta_t^i$,
with a total number of $(32+bd)$ bits. It is dramatically smaller than $32d$ bits required by unquantized gradient.
Based on the two data-parallel frameworks introduced in Section \ref{sec:data-parallel-ar}, there are three possible communication schemes: 
\\
\\
\noindent
\textbf{(a.)} (\textbf{Broadcast})
For each worker $i$, the scale factor $\delta_t^i$ is computed using (\ref{eq:delta_t_i}).
%
%
When receiving low-precision gradients from other workers, it first recovers their full-precision representations and calculates
$\tilde{u}_t$ using equation~(\ref{eq:tilde_u_t}).
 In this case, precision loss only happens in calculating $\tilde{u}_{t}^i$, and the resulting communication cost is
 $(32+bd)N(N-1)$ bits.
\begin{equation}
  \label{eq:tilde_u_t}
  \tilde{u}_t = \frac{1}{N}\sum_{i=1}^N \tilde{u}_{t}^i.
\end{equation}
\\
\noindent
\textbf{(b.)} (\textbf{Parameter-server without re-quantizing})
%
The adding operation can only be conducted among low-precision numbers with the same scale factor.
Therefore,
$2N$ floating numbers need to be transmitted between server and workers to unify the same scale factor $\delta_t = \max_i\{\delta_t^i\}$.
Then all workers adopt $\delta_t$ in Step $7$, i.e., $\delta_t^i = \delta_t$ for each $i$.
When server calculates the average gradients $\tilde{u}_t$ using (\ref{eq:tilde_u_t}), it allows for an overflow of $\lceil\log{N}\rceil$ bits, and then send $\tilde{u}_t$ to workers.
The total communication overhead is $N(64+2bd+d\lceil\log{N}\rceil)$ bits,
and the precision loss also comes from calculating  $\tilde{u}_{t}^i$.
\\
\\
\noindent
\textbf{(c.)} (\textbf{Parameter-server with re-quantizing})
Communication scheme is the same as $\textbf{(b)}$,
except for
\begin{equation}
  \label{eq:re-quan}
   \tilde{u}_t = Q_{(\delta_t,b)} (\frac{1}{N}\sum\limits_{i=1}^N\tilde{u}_t^i), \quad \delta_t = \max_i\{\delta_t^i\},
\end{equation}
i.e., server re-quantizes $(\frac{1}{N}\sum_i\tilde{u}_t^i)$ to a $b$-bit low-precision representation with scale factor $\delta_t$ (to prevent bit overflow).
Note that the quantization in (\ref{eq:re-quan}) is unbiased since all $u_t^i$\textprime s are also quantized with $\delta_t$.
In this case, there exist two levels of quantization errors and the total communication complexity equals to $N(64+2bd)$.
The comparisons of three communication schemes are summarized in Table~\ref{table:com-scheme}.

\subsection{Theoretical Analysis}
\label{sec:LPC-SVRG-ana}
Based on the above low-precision representation and communication schemes,
in this section, we provide the convergence analysis of Algorithm~\ref{Algm:LPC-SVRG}.
We begin with several lemmas which bound the variance of low-precision gradient. 
\begin{lemma} If Assumptions~\ref{assum:unbia}, \ref{assum:lip} hold,\\
  \label{lemma:biased-variance}
  (i).
  under communication scheme \textbf{(a)} or \textbf{(b)} with $\tilde{u}_t = \frac{1}{N}\sum_{i} \tilde{u}_{t}^i$, we have
  \begin{equation*}
    \begin{aligned}
  \mathbf{E}||v_t^{s+1} - \nabla f(x_t^{s+1})||^2 \leq
   2L^2 \Big[ \frac{d\lambda^2} {4(2^{b-1}-1)^2} + d_\lambda(1-\lambda)^2 + \frac{1}{NB} \Big] \mathbf{E} ||x_t^{s+1} - \tilde{x}^s||^2;
\end{aligned}
\end{equation*}
  (ii). under communication scheme \textbf{(c)} with
$\tilde{u}_t = Q_{(\delta_t,b)} (\frac{1}{N}\sum\limits_{i=1}^N\tilde{u}_t^i)$,
  we obtain a variance bound of
\begin{equation*}
  \begin{aligned}
    \label{eq:re-quan-variance}
\mathbf{E}||v_t^{s+1} - \nabla f(x_t^{s+1})||^2 \leq
2L^2\Big[ \frac{3d\lambda^2} {8(2^{b-1}-1)^2} + d_\lambda(1-\lambda)^2 + \frac{1}{NB} \Big] \mathbf{E} ||x_t^{s+1} - \tilde{x}^s||^2,
\end{aligned}
\end{equation*}
\noindent
where $d_{\lambda} = \max_i\{d_{\lambda}^i\}$.
$d_{\lambda}^i$ is the number of coordinates in $u_t^i$ exceeding $\dom(\delta_t^i,b)$.
Note that $\delta_t^i = \delta_t$ in schemes \textbf{(b)} and \textbf{(c)}.
\end{lemma}
\noindent
\emph{Remarks.} The coefficients of variance bound in Lemma~\ref{lemma:biased-variance}, e.g., \emph{(i)}
can be decomposed into three items, i.e.,
$A = \frac{d\lambda^2} {4(2^{b-1}-1)^2}$, $B = d_\lambda(1-\lambda)^2$, $C = \frac{1}{NB}$,
where $A$ and $B$ come from quantization and the last item $C$ inherits the gradient variance of SVRG \cite{Reddi2016Fast, johnson2013accelerating}.
If without quantization, $A=B=0$, and
Lemma~\ref{lemma:biased-variance} becomes equivalent to Lemma $3$ in \cite{Reddi2016Fast}.
%
%
Moreover, Lemma~\ref{lemma:biased-variance} theoretically shows the benefit of quantization with gradient clipping, i.e.,
with proper choice of $\lambda < 1$,  $(A|_{\lambda <1}+B)<A|_{\lambda=1}$.
Such condition can always hold in the case where gradients follow a Gaussian distribution.
Specifically, when $A|_{\lambda=1}$ is dominating, a smaller $\lambda$ can effectively shrink its value while fewer coordinates, i.e., very small $d_{\lambda}$, comes into $B$.
\begin{theorem}
  \label{thm:LPC-SVRG}
  Suppose Assumptions~\ref{assum:unbia}, \ref{assum:lip} hold, $h(x)$ is convex, $T=Sm$, $\eta=\frac{\rho}{L}$, $\rho < \frac{1}{2}$,
  and all parameters satisfy
  \begin{equation}
    \label{eq:con-1}
    8m^2\rho^2\Big[ \frac{d\lambda^2}{4(2^{b-1}-1)^2} + d_{\lambda}(1-\lambda)^2  +\frac{1}{NB} \Big] + \rho \leq 1
  \end{equation}
  when using communication scheme \textbf{(a)} or \textbf{(b)} ($d_\lambda$ is defined in Lemma~\ref{lemma:biased-variance}).
 For the output $x_{out}$ of Algorithm~\ref{Algm:LPC-SVRG}, we have
  \begin{equation}
    \label{eq:LPC-SVRG}
    \mathbf{E}||G_\eta(x_{out})||^2 \leq \frac{2L(P(x^0)-P(x^*))}{\rho(1-2\rho)T}.
  \end{equation}
Moreover, if communication scheme \textbf{(c)} is adopted, the same rate (\ref{eq:LPC-SVRG}) can be obtained if
  \begin{equation}
      \label{eq:con-2}
  8m^2\rho^2\Big[ \frac{3d\lambda^2}{8(2^{b-1}-1)^2} + d_{\lambda}(1-\lambda)^2 +\frac{1}{NB} \Big] + \rho \leq 1.
  \end{equation}
\end{theorem}
\begin{corollary}
  \label{coro:LPC-SVRG}
  If $b=O(\log{\sqrt{d}})$ and $\eta=O(1/mL)$, there exists $\lambda \in (0,1]$ that guarantees constraints (\ref{eq:con-1}) and (\ref{eq:con-2}).
\end{corollary}
%
%
%
%
%
\noindent
The above theorem shows the same asymptotic convergence rate (as Theorem~\ref{thm:full-precision-svrg}) can be achieved with low-precision representation.
In the following section, we propose a faster and communication-efficient algorithm with \emph{double sampling}.

\begin{algorithm}[tb]
   \caption{Accelerated LPC-SVRG: \textbf{ALPC-SVRG}  (for each worker $i$, strongly convex case)  }
\begin{algorithmic}[1]   \label{Algm:ALPC-SVRG-s}
   \STATE {\bfseries Input:} $S$, $m$, $\lambda$, $B$, $\tau_1$, $\tau_2$, $\alpha$, $y_0 = z_0 = x_0 = \tilde{x}^0$;
   \FOR{$s=0,1,...,S-1$}
       \STATE \textbf{Communication step $1$:} cooperates with other workers to compute $\nabla f(\tilde{x}^s)$;
       \FOR{$t=0$ {\bfseries to} $m-1$}
           \STATE  $k \gets (sm+t)$;
           \STATE  $x_{k+1} = \tau_1 z_k + \tau_2 \tilde{x}^s + (1-\tau_1 - \tau_2) y_k$;
           \STATE  uniformly and independently samples (with replacement) a mini-batch $i_B$ with size B to calculate\\
            $u_{k+1}^i = \frac{1}{B} \sum\limits_{a \in i_B} \big[ \nabla f_a(x_{k+1}) - \nabla f_a(\tilde{x}^s) \big]$;
           \STATE
           \textbf{Quantization step:} $\tilde{u}_{k+1}^i = Q_{(\delta_{k+1}^i, b)} (u_{k+1}^i)$;
           \STATE  
           \textbf{Communication step $2$:} cooperates with other workers to compute $\tilde{u}_{k+1} = \frac{1}{N}\sum\limits_{i=1}^N \tilde{u}_{k+1}^i$;
           \STATE $v_{k+1} = \tilde{u}_{k+1} + \nabla f(\tilde{x}^s)$;
           \STATE
           \textbf{Double sampling:} independently (to $i_B$) samples a mini-batch $J$ with size $B$ to compute
           $\hat{v}_{k+1} = \frac{1}{B}\sum\limits_{j \in J} \Big(\nabla f_j(x_{k+1}) - \nabla f_j(\tilde{x}^s) \Big) +\nabla f(\tilde{x}^s)$;
       \STATE $ y_{k+1} = \argmin_{y} \{ 2L||y-x_{k+1}||^2 + \langle v_{k+1},y \rangle + h(y) \}$;\quad $\diamondsuit$ \emph{with low-precision gradient}
      \STATE $ z_{k+1} = \argmin_{z}\{ \frac{1}{2\alpha} ||z-z_k||^2 + \langle \hat{v}_{k+1}, z \rangle + h(z)\}$; \quad $\diamondsuit$ \emph{with local full-precision gradient}
       \ENDFOR
      \STATE $\tilde{x}^{s+1} = \Big(\sum\limits_{t=0}^{m-1}(1+\alpha\sigma)^t\Big)^{-1} \cdot \Big(\sum\limits_{t=0}^{m-1}(1+\alpha\sigma)^t \cdot y_{sm+t+1}\Big)$;
  \ENDFOR
   \STATE {\bfseries Output:} $\tilde{x}^{S}$.
\end{algorithmic}
\end{algorithm}
\section{ACCELERATED LOW-PRECISION ALGORITHM}
%
Recently momentum or Nesterov \cite{nesterov1983method} technique has been successfully combined with SVRG to achieve faster convergence in real-world applications~\cite{allen2017katyusha,hien2016accelerated,lin2015universal}.
%
%
In this section, we propose an accelerated method named ALPC-SVRG, based on Katyusha momentum~\cite{allen2017katyusha}, to obtain both faster running speed and high communication-efficiency.
As shown in Algorithm~\ref{Algm:ALPC-SVRG-s}, at the beginning of inner iterations,
we initialize $x_{k+1}$ as a convex combination of the reference point $\tilde{x}^s$ and two auxiliary variables $z_k$ and $y_k$.
We propose \emph{double sampling} to compute the stochastic gradients for $z_k$ and $y_k$,
where the gradient $v_{k+1}$ used in $y_k$-update is the same as that of Algorithm~\ref{Algm:LPC-SVRG},
and a new gradient $\hat{v}_{k+1}$ with independent \emph{double sampling} is applied in computing $z_{k+1}$.
Note that $\hat{v}_{k+1}$ has full-precision and is calculated locally in each worker without communication.
To make sure the consistency of model $x$ in all workers, we assume the mini-batches $J$ used by
all workers are the same.
This can be achieved by setting the identical random seeds for $J$. 
\\
\\
\noindent
In this section, quantization with gradient clipping is also considered and $\lambda$ is the clipping parameter, i.e., $\delta_{k+1}^i = \frac{\lambda||u_{k+1}^i||_{\infty}}{2^{b-1}-1}$.
%
%
Without loss of generality, we use the communication scheme \textbf{(a)}  for the analysis below,
and it is easy to extend the following conclusions to the other two communication schemes.
\begin{theorem}
  \label{thm:ALPC-SVRG-s}
  Suppose Assumptions \ref{assum:unbia}, \ref{assum:lip} hold and each $f_i(x)$ is convex, $h(x)$ is $\sigma$-strongly convex,
  i.e., there exists $\sigma > 0$ such that $\forall x, y$,
  \begin{equation*}
    h(y) \geq h(x) + \langle\nabla h(x), y-x\rangle + \sigma/2||y-x||^2.
  \end{equation*}
  Denote $\zeta \triangleq \frac{d\lambda^2}{4(2^{b-1}-1)^2} + d_{\lambda}(1-\lambda)^2 + \frac{1}{NB}$,
  where
  $d_{\lambda} = \max_i\{d_{\lambda}^i\}$,
  $d_{\lambda}^i$ is the number of coordinates in $u_{k+1}^i$ exceeding $\dom(\delta_{k+1}^i,b)$.
  Let $\alpha = \frac{1}{6\tau_1 L}$, $m \leq \frac{3L}{2\sigma}$, if $\tau_1$, $\tau_2$ satisfy
  \begin{equation*}
    \tau_1 = \sqrt{\frac{m\sigma}{6L}}, \quad
    \tau_2 = \frac{5\zeta}{3} + \frac{1}{2B} \leq \frac{1}{2}.
  \end{equation*}
%
Then under Algorithm~\ref{Algm:ALPC-SVRG-s}, we obtain
%
  \begin{equation*}
    \begin{aligned}
   \mathbf{E}\Big[ P(\tilde{x}^S) &- P(x^*)\Big]
   \leq  O\Big( \big(1+\sqrt{ \frac{\sigma}{Lm} }\big)^{-Sm}\Big) \Big(P(x_0) - P(x^*)\Big).
\end{aligned}
  \end{equation*}
\end{theorem}
%
\noindent
It can be verified that
if $b=O(\log{\sqrt{d}})$, there exist $\lambda \in (0,1]$ and mini-batch size $B$ that make sure $\tau_2 \leq  \frac{1}{2}$,
and we obtain the same convergence rate as the full-precision algorithm (Katyusha, Case 1) in \cite{allen2017katyusha}.
\\
\\
\noindent
The following theorem provides the convergence rate for ALPC-SVRG without assuming strong convexity.
%
%
%
%
As shown in Algorithm~\ref{Algm:ALPC-SVRG-g}, we update the value of $\tau_1$ and $\alpha$ at each epoch $s$, and $\tau_2$ is the same as Algorithm~\ref{Algm:ALPC-SVRG-s}.
The update of $\tilde{x}^{s+1}$ is also adjusted to conduct telescoping summation.
Based on the same communication scheme and low-precision representations,  
we have the following result.

\begin{algorithm}[tb]
   \caption{Accelerated LPC-SVRG: \textbf{ALPC-SVRG} (for each worker $i$, general convex case) }
\begin{algorithmic}[1]   \label{Algm:ALPC-SVRG-g}
   \STATE {\bfseries Input:} $S$, $m$, $\lambda$, $B$, $\tau_2$, $y_0 = z_0 = x_0 = \tilde{x}^0$;
   \FOR{$s=0,1,...,S-1$}
       \STATE updates $\tau_{1,s} = \frac{2}{s+4}$, $\alpha_s = \frac{1}{6L\tau_{1,s}}$;
       \STATE performs the same Steps $3$-$14$ as in Algorithm~\ref{Algm:ALPC-SVRG-s};
      \STATE $\tilde{x}^{s+1} = \frac{1}{m} \sum\limits_{t=1}^{m}y_{sm+t}$;
  \ENDFOR
   \STATE {\bfseries Output:} $\tilde{x}^{S}$.
\end{algorithmic}
\end{algorithm}

\begin{theorem}
Denote $\zeta \triangleq \frac{d\lambda^2}{4(2^{b-1}-1)^2} + d_{\lambda}(1-\lambda)^2 + \frac{1}{NB}$.
Let $\tau_2 = \frac{5\zeta}{3} + \frac{1}{2B}$, if Assumptions~\ref{assum:unbia}, \ref{assum:lip} hold, each $f_i(x)$ and $h(x)$ are convex and $\tau_2 \leq \frac{1}{2}$,
  then under Algorithm~\ref{Algm:ALPC-SVRG-g}, we have
  \begin{equation*}
  \begin{aligned}
  \mathbf{E} \Big[ P(\tilde{x}^S) - P(x^*)\Big]
  \leq O(\frac{1}{mS^2}) \Big[ m(P(x_0) - P(x^*)) + L ||x_0 - x^*||^2 \Big].
  \end{aligned}
\end{equation*}
\end{theorem}

\begin{figure*}[!ht]
\begin{center}
 \subfigure[Syn-512]{
\minipage{0.29\textwidth}
  \includegraphics[width=1.1\linewidth]{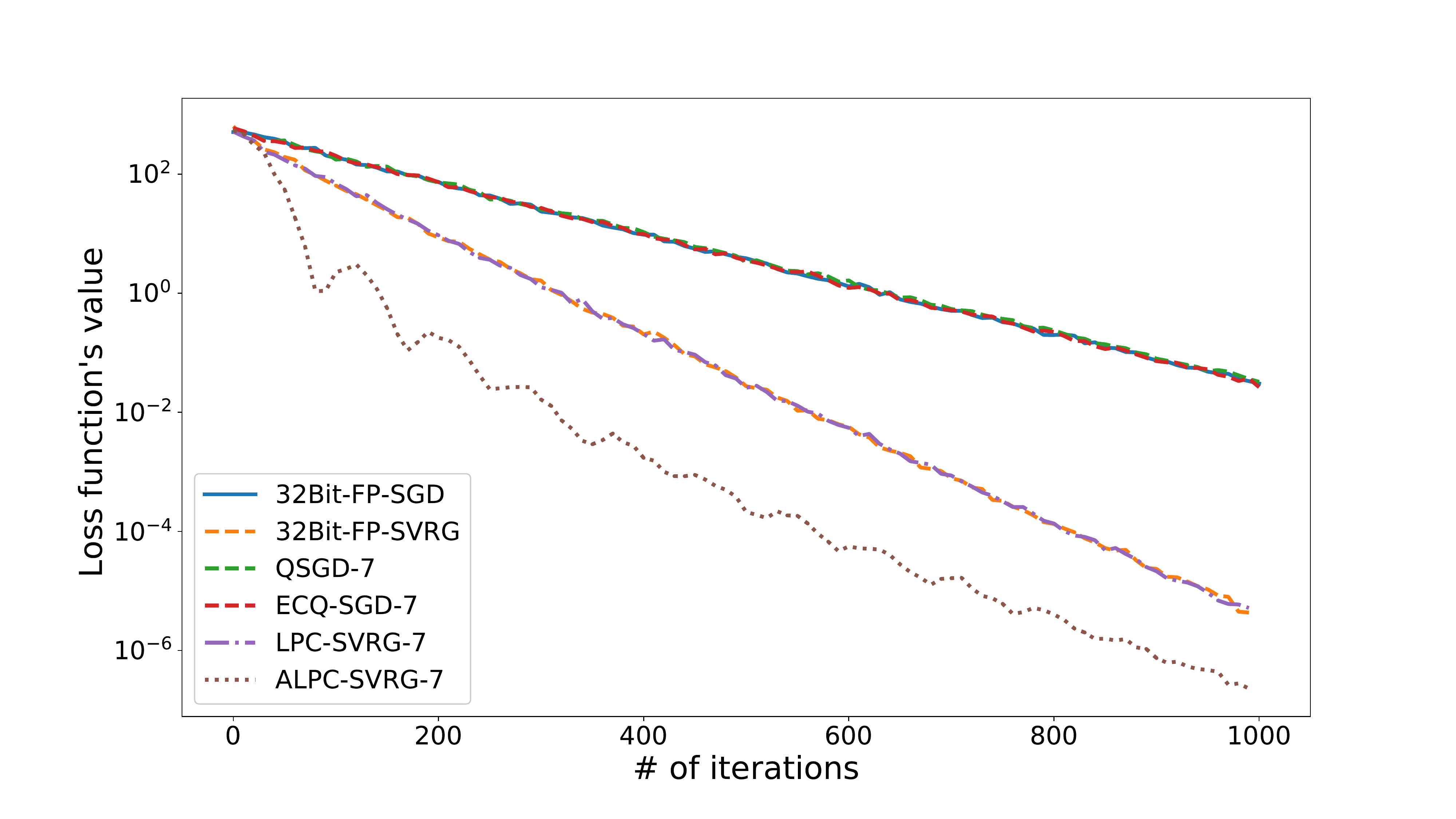}
 \label{fig:512}
\endminipage
\hfill
}
 \subfigure[Syn-1024]{
\minipage{0.29\textwidth}
  \includegraphics[width=1.1\linewidth]{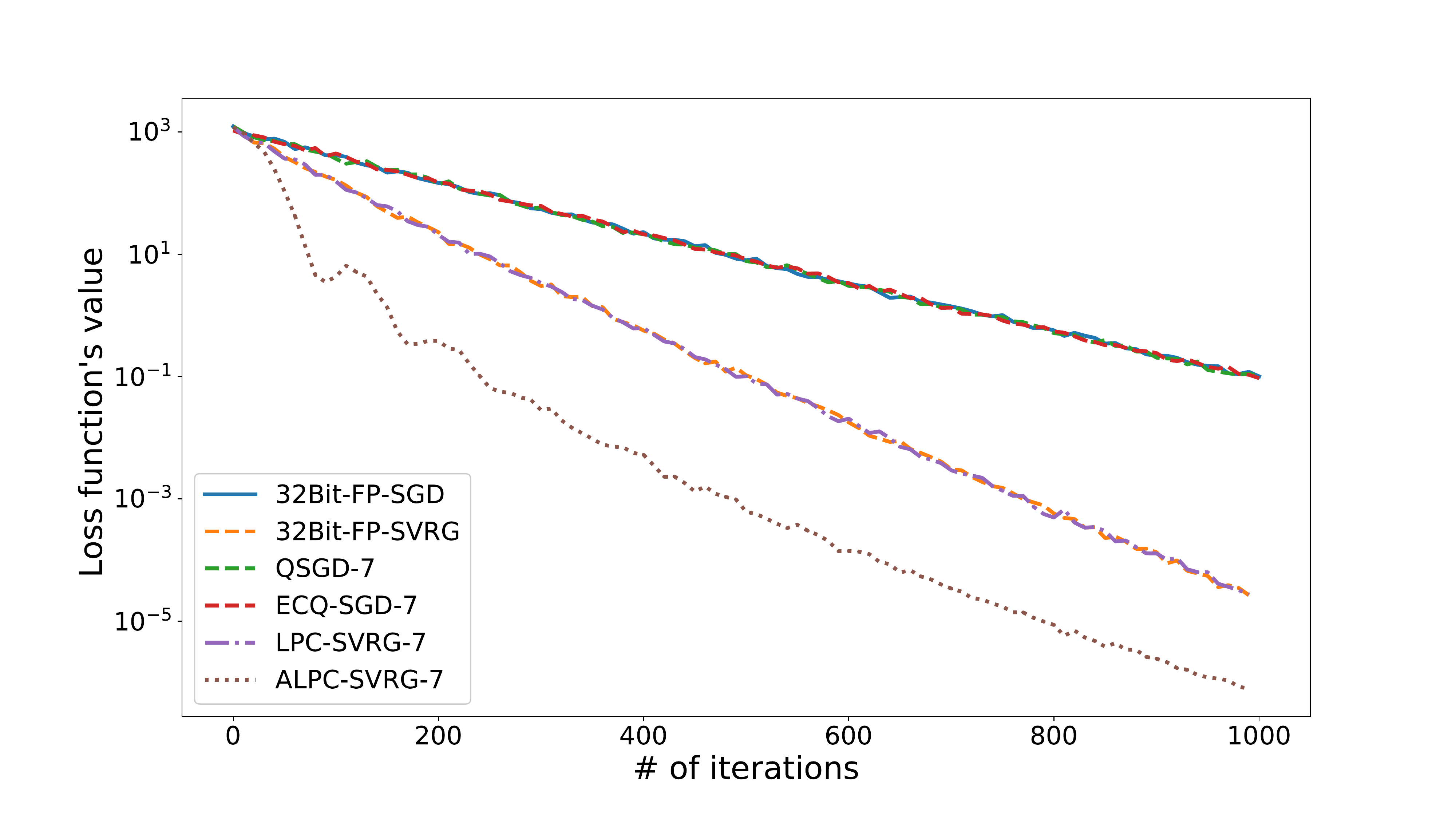}
  \label{fig:1024}
\endminipage
\hfill
}
\subfigure[Syn-20k]{
\minipage{0.32\textwidth}
 \includegraphics[width=\linewidth]{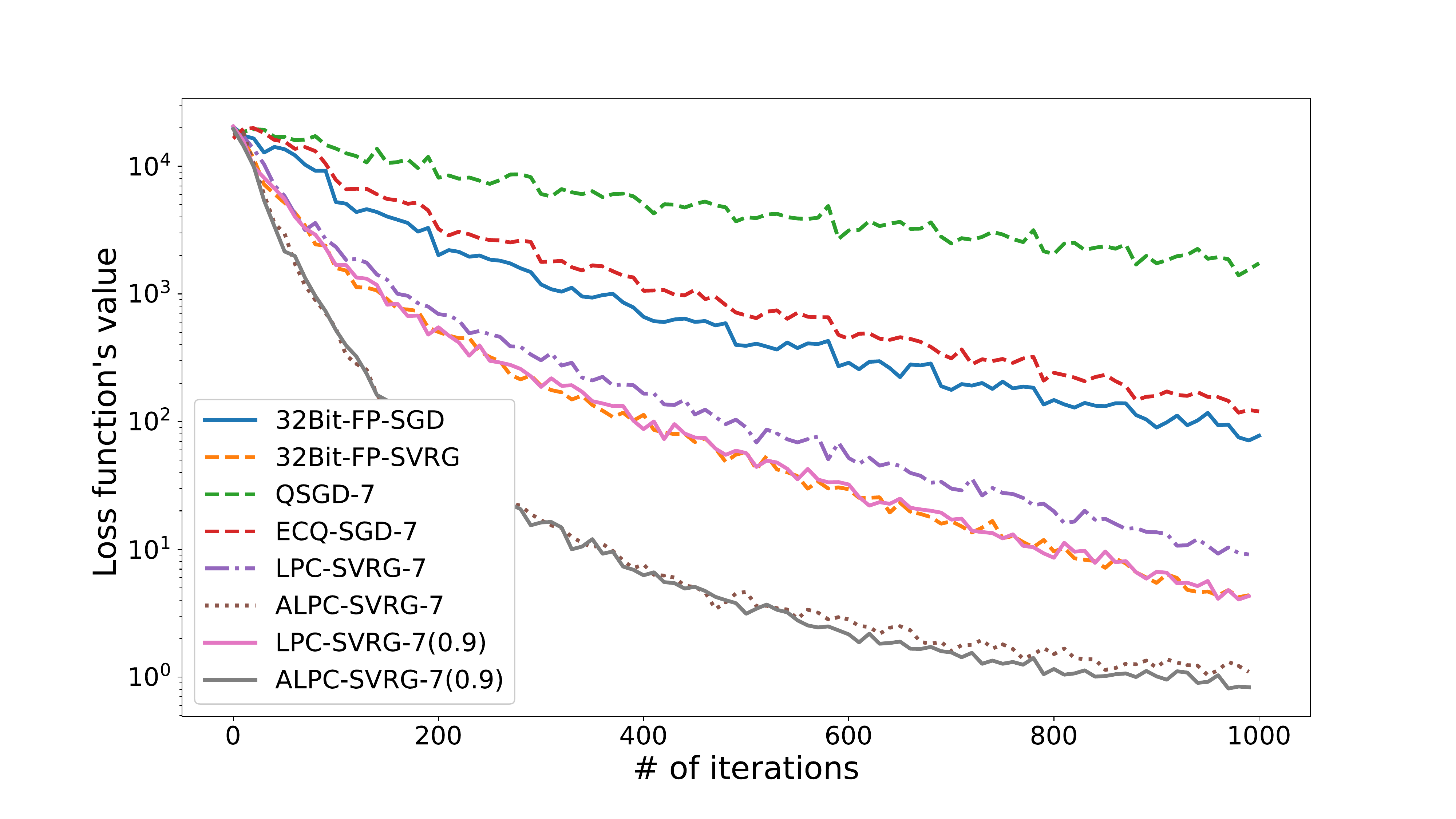}
 \label{fig:20k}
\endminipage
\hfill
}
\caption{The training curves of compared algorithms on three synthetic datasets.
Here the suffix denotes the number of quantization levels and the number in bracket is the value of clipping parameter.}
\label{fig:training-loss-512-1024}
\end{center}
\end{figure*}

\section{EXPERIMENTS}
In this section, we conduct extensive experiments to validate the effectiveness of our
proposed
algorithms.
Firstly, we evaluate linear models on various datasets
and then extend our algorithms to train deep neural networks on datasets CIFAR-10 \cite{krizhevsky2009learning} and ILSVRC-12 \cite{russakovsky2015imagenet}.
In all evaluations, the communication scheme $\textbf{(a)}$ described in Section~\ref{sec:Communication Schemes and Complexities} is adopted.

\subsection{Evaluation on Linear Model}
\label{sec:linear-model}
We begin with linear regression on three synthetic datasets: Syn-512, Syn-1024, Syn-20k\footnote{These datasets are generated using the same method as in \cite{wu2018error} with \emph{i.i.d} random noise.},
each containing 10k, 10k, 50k training samples.
Here the suffix denotes the dimension of model parameters.
We compare LPC-SVRG and ALPC-SVRG with 32 bits full-precision algorithms:
 SGD, SVRG~\cite{johnson2013accelerating},
and the state-of-the-art low-precision algorithms:
QSGD~\cite{alistarh2017qsgd},
ECQ-SGD~\cite{wu2018error}.
Among the above algorithms, QSGD and ECQ-SGD are the low-precision variants of SGD,
and LPC/ALPC-SVRG is based on full-precision SVRG~\cite{johnson2013accelerating,Reddi2016Fast}.
%
%
Similar to related works \cite{zhang2017zipml,allen2017katyusha,Reddi2016Fast},
we use diminishing/constant learning rate (\emph{lrn\_rate} for short) for SGD/SVRG based algorithms respectively.
In our experiments, we have tuned the \emph{lrn\_rates} for the best performance of each algorithm.
The hyper-parameters in ALPC-SVRG are setting to be $\tau_1 = \frac{2}{s+4}$, $\alpha = \emph{lrn\_rate}/\tau_1$, $\tau_2 = \frac{1}{2}$.
\\
\\
\noindent
For low-precision algorithms, we use max norm in the scaling factor
and  adopt the entropy encoding scheme \cite{huffman1952method} to further reduce the communication overhead. 
Here for the consistency of notation, we adopt \emph{levels} mentioned in \cite{alistarh2017qsgd,wu2018error} to represent the extent of low-precision representation.
As mentioned in Section~\ref{sec:Low-Precision with Gradient Clipping}, the quantization \emph{levels} equal to the number of positive points in $\dom(\delta,b)$.
\\
\\
\noindent
\textbf{Convergence.}
Figure~\ref{fig:512} and~\ref{fig:1024} depict the training curves of six algorithms on datasets syn-512 and syn-1024,
where the suffix of low-precision algorithm denotes the number of quantization levels.
It shows that all low-precision algorithms have comparable convergence to their corresponding full-precision algorithms,
which verifies the redundancy of 32-bit representation.
On the other hand, we can see that LPC-SVRG converges faster than QGD and ECQ-SGD with the same quantization levels.
And the optimal convergence of ALPC-SVRG is validated in these two figures.
\\
\\
\noindent
\textbf{Gradient clipping.}
We conduct gradient clipping on a larger dataset syn-20k,
where for LPC/ALPC-SVRG-$7(0.9)$, $7$ represents the quantization levels and $0.9$ is the clipping parameter, i.e., $\lambda$.
The value of $\lambda$ is determined by cross validation and we only conduct it once.
As shown in Figure~\ref{fig:20k},
LPC-SVRG without clipping converges slower than SVRG because the precision loss becomes significant for the larger dataset syn-20k.
The same circumstance is also presented in QSGD and ECQ-SGD.
In this case, gradient clipping comes into effect.
With more dense quantization points, it can reduce the quantization error and achieve better performance.
\begin{figure}[h]
  \centering
 \includegraphics[width=0.5\linewidth]{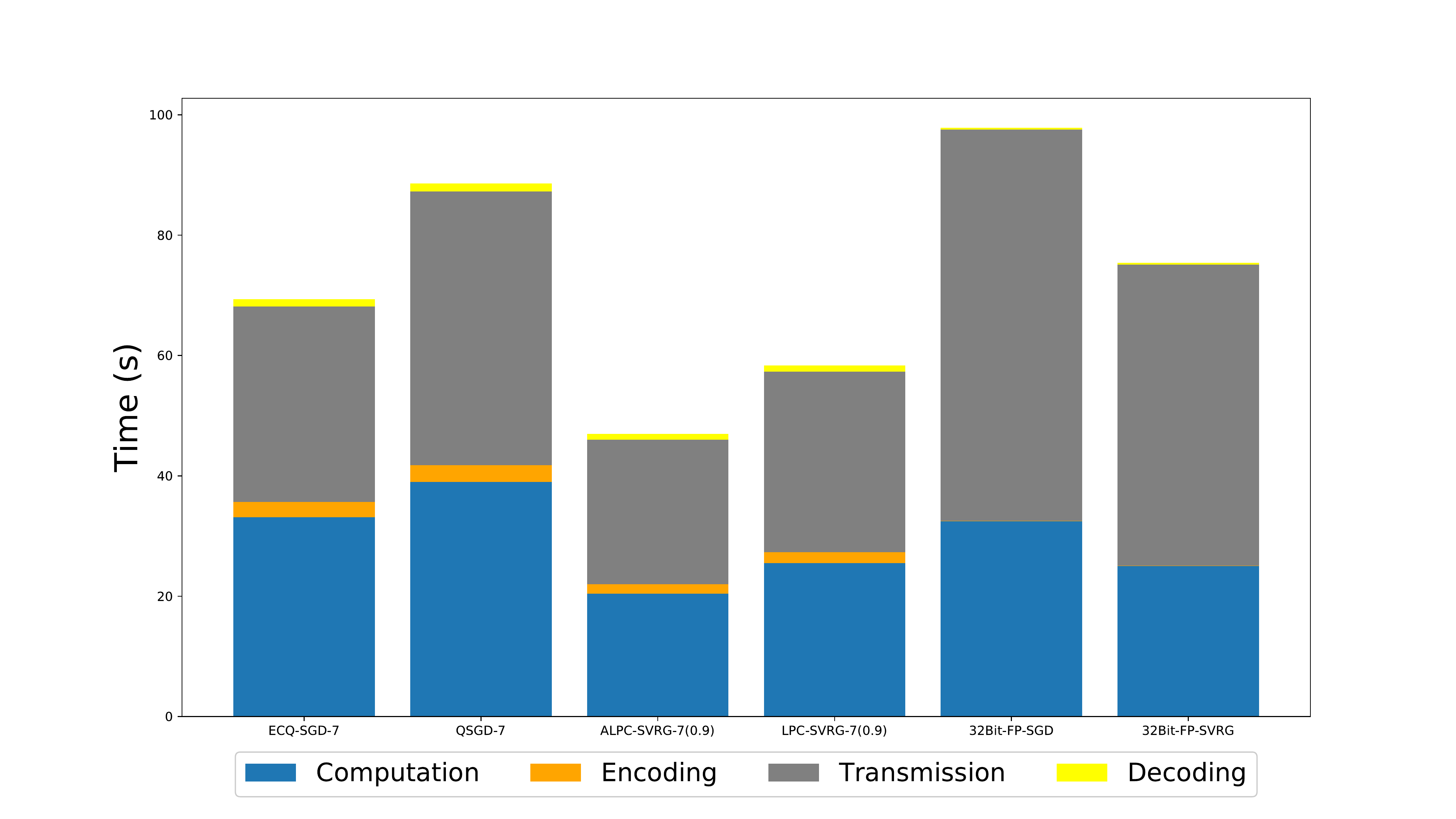}
\caption{Comparisons on the time consumption of different algorithms on dataset syn-20k.}
\label{fig:time}
\end{figure}
\\
\\
\noindent
\textbf{Communication time.}
To better demonstrate the effectiveness of LPC/ALPC-SVRG, we analyze the time consumption for different algorithms on dataset syn-20k.
Specifically, we decompose the total running time into 4 parts, including \emph{computation}, \emph{encoding}, \emph{transmission}, \emph{decoding}, and the time is reported till similar convergence (i.e. first training loss below $10^2$).
As shown in Figure~\ref{fig:time},
compared to full precision SVRG,
LPC-SVRG saves the total running time by significantly reducing the transmission overhead,
and ALPC-SVRG requires the fewest computation and transmission time among all algorithms.
\\
\\
%
\begin{figure}
\begin{minipage}{0.5\linewidth}
  \centering
  \includegraphics[width=2in]{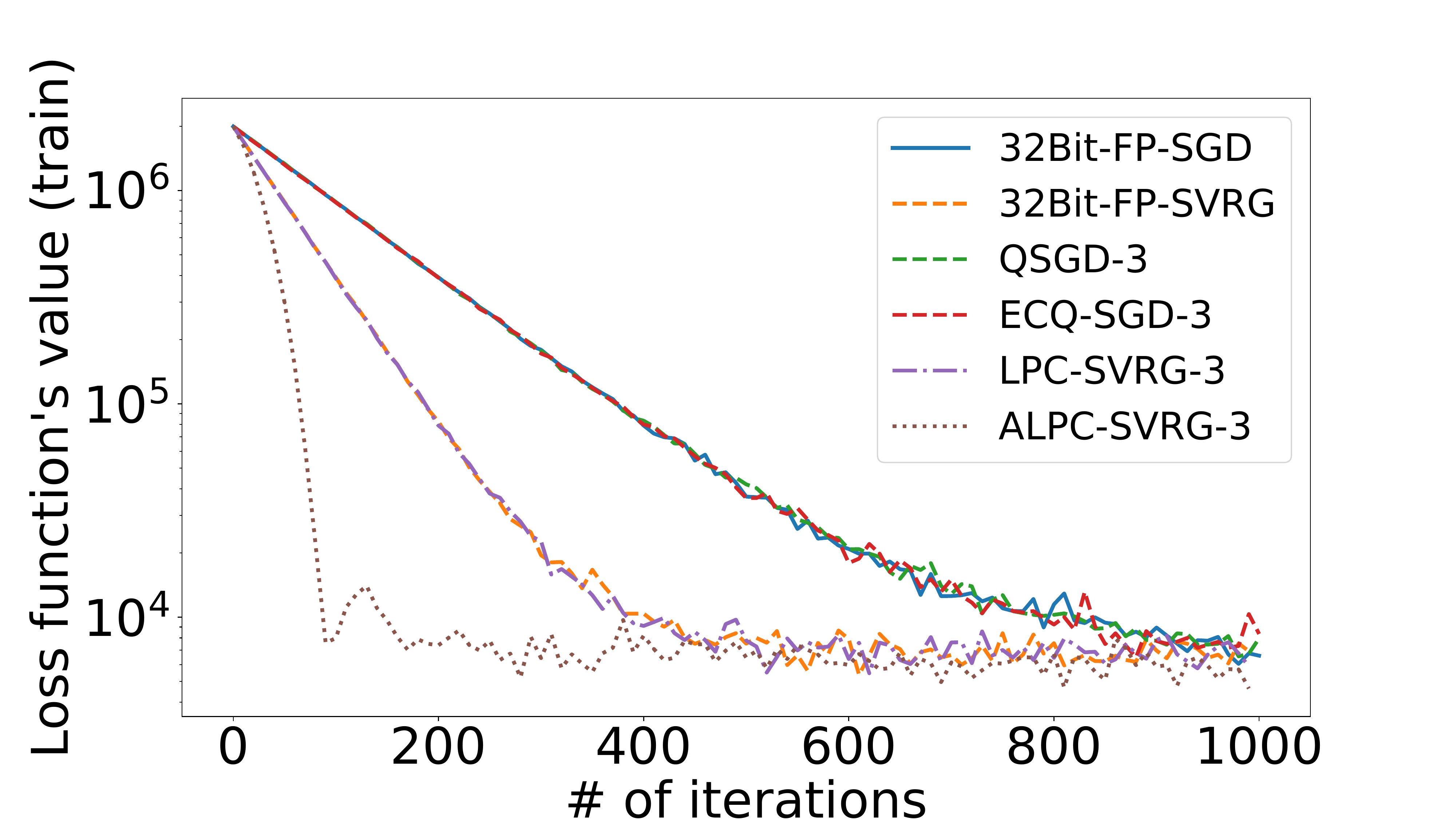}
\end{minipage}%
\begin{minipage}{0.5\linewidth}
  \centering
  \includegraphics[width=2in]{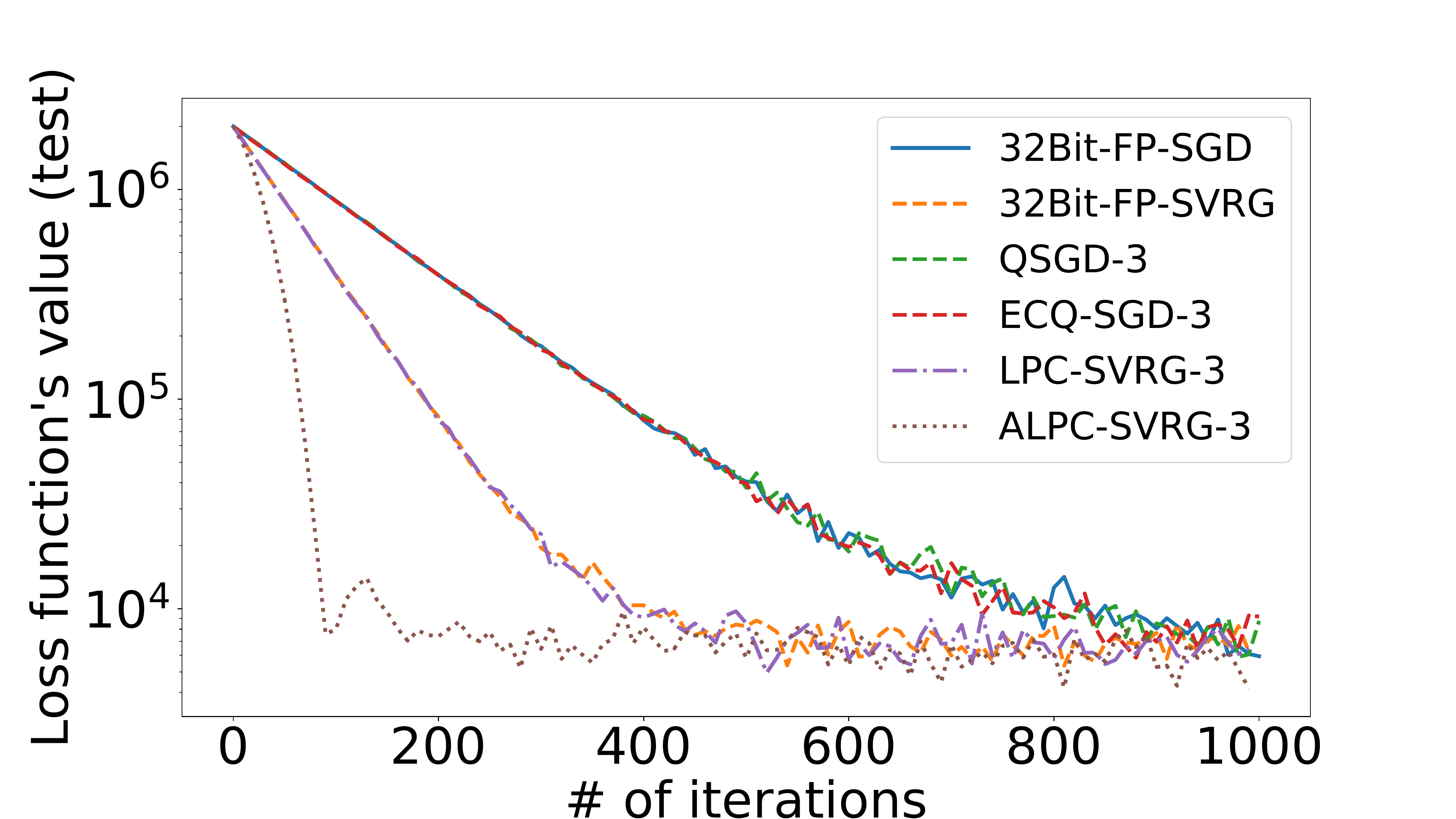}
\end{minipage}
\caption{Validation on dataset \emph{YearPredictionMSD} (left: training loss, right: test loss).}
\label{fig:year}
\end{figure}
%
%
\begin{table}[h]
\caption{Communication costs of different algorithms for linear regression on dataset \emph{YearPredictionMSD}.
All low-precision algorithms below use three quantization levels.}
\label{table:comm-cost}
\begin{center}
\begin{tabular}{lccc}
\toprule
\textbf{Algorithm}  &\textbf{Training loss}  &\textbf{\# of bits} &\textbf{Ratio}\\
\midrule
32-bit-FP-SGD         & 5.74$e3$ & 3.25$e6$ & $-$\\
32-bit-FP-SVRG        & 5.55$e3$ & 1.41$e6$ & $2.3\times$\\
QSGD                  & 5.96$e3$ & 1.61$e5$ & $20.19\times$\\
ECQ-SGD               & 5.90$e3$ & 1.60$e5$ & $20.31\times$\\
LPC-SVRG              & 5.80$e3$ & 7.04$e4$ & $46.16\times$ \\
ALPC-SVRG             & 5.86$e3$ & 3.50$e4$ & $92.86\times$\\
\bottomrule
\end{tabular}
\end{center}
\end{table}
\noindent
\textbf{Real-world dataset and communication overhead.}
We further evaluate linear regression on public dataset \emph{YearPredictionMSD} \cite{chang2011libsvm}.
In this case, fewer (three) levels  are used since the feature dimension is smaller (i.e., 90).
As shown in Figure~\ref{fig:year}, all low-precision algorithms achieve similar convergence as their corresponding full-precision algorithms.
Moreover, the empirical results validate the fast convergence of LPC/ALPC-SVRG.
\\
\\
\noindent
In Table~\ref{table:comm-cost}, we record the total number of communication costs of various algorithms for achieving similar training loss.
It shows LPC-SVRG and ALPC-SVRG can save up to $46.16\times$ and $92.86\times$ communication costs compared with the benchmark,
which validates the effectiveness of our proposed algorithms.

\subsection{Evaluation on Deep Learning Model}
We further extend our algorithms to train deep neural networks.
The public dataset CIFAR-10 \cite{krizhevsky2009learning} is used, with 50k training and 10k test images.
We setup experiments on TensorFlow \cite{abadi2016tensorflow} with ResNet-20 model \cite{he2016identity}, and
%
 adopt a decreasing learning rate, i.e.,
starting from $0.1$ and divided by 10 at 40k and 60k iterations.
All low-precision algorithms use the same quantization levels, i.e., 7.
%
The other training hyper-parameters are tuned to achieve similar loss as the benchmark full-precision SGD.
\begin{figure}
\begin{minipage}{0.5\linewidth}
  \centering
  \includegraphics[width=2.3in]{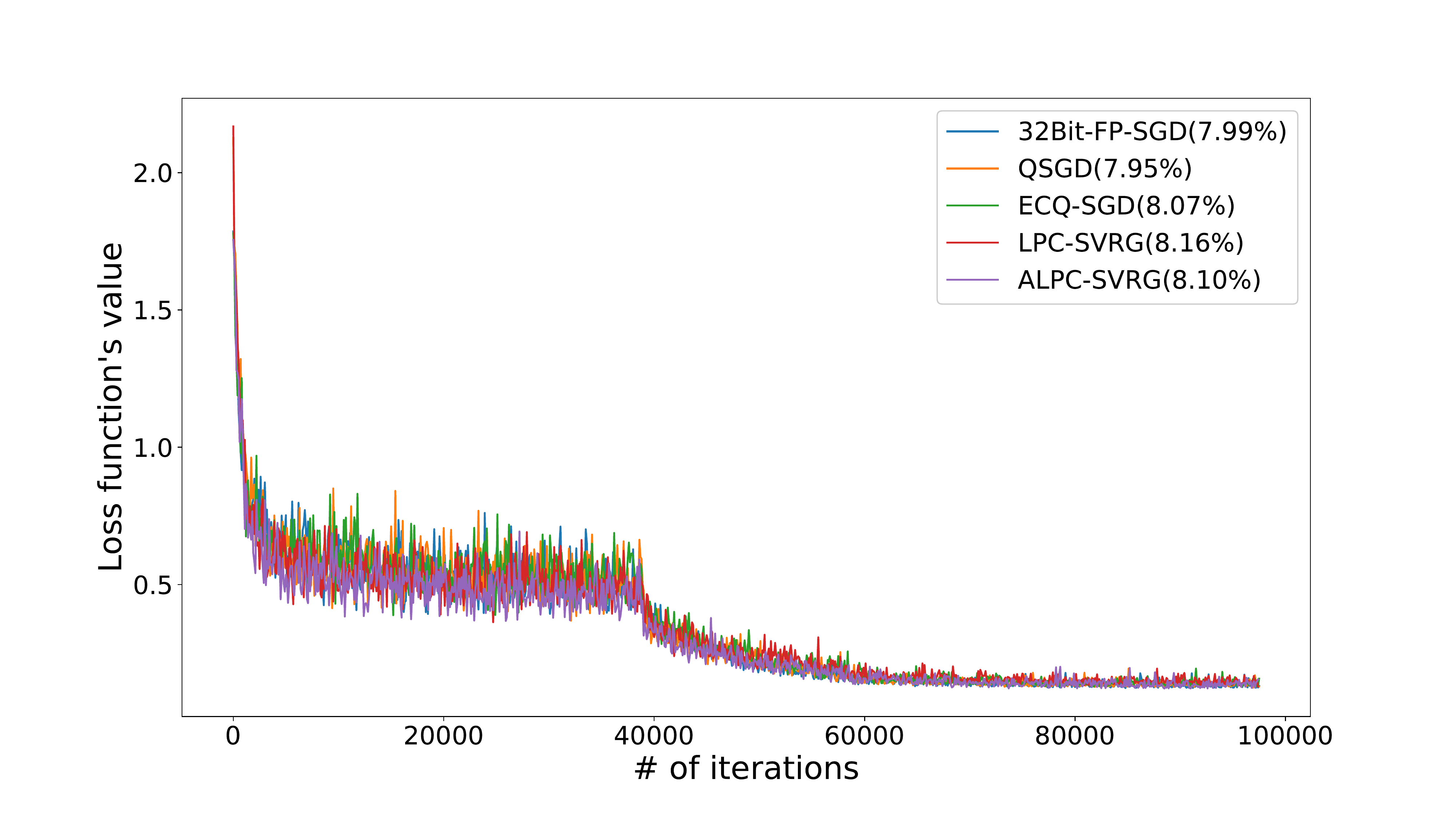}
\end{minipage}%
\begin{minipage}{0.5\linewidth}
  \centering
  \includegraphics[width=2.3in]{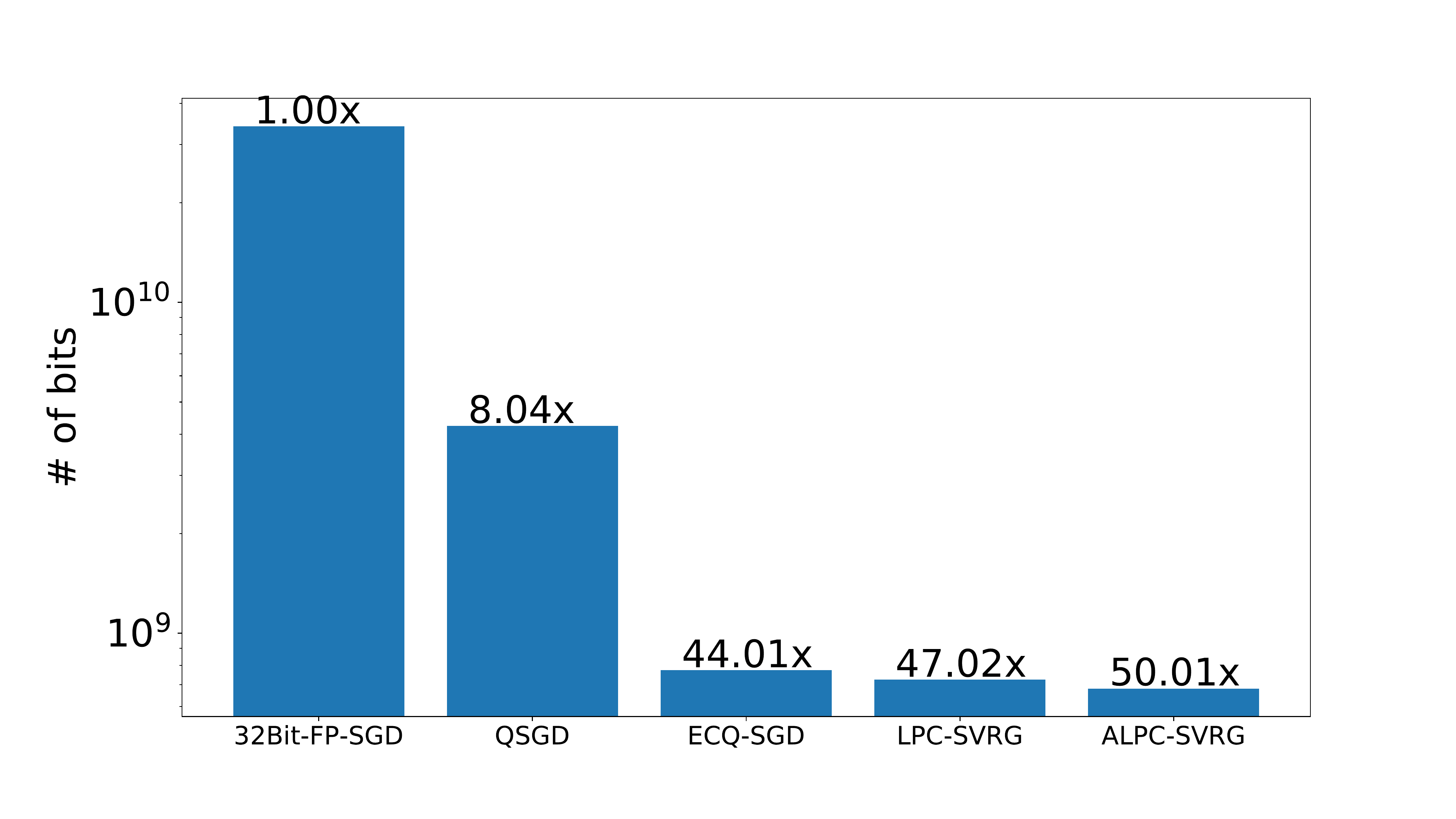}
\end{minipage}
\caption{Validation on CIFAR-10 with ResNet-20 model,
left: training curves, right: the amount of transmitted bits.}
\label{fig:resnet}
\end{figure}


\noindent
Figure~\ref{fig:resnet} shows the training statistics of each algorithm.
The left graph reports the training curves versus iterations and
the values in bracket are test classification error.
These curves show that our algorithms converge with a faster speed.
The right graph plots the total amount of transmitted bits of compared algorithms (till similar convergence, i.e., training loss first below 0.17).
In this experiment, a clipping parameter with value 0.85 is adopted in LPC-SVRG and ALPC-SVRG.
We discover that our proposed algorithms effectively save the communication overhead.
\begin{figure}[h]
  \centering
 \includegraphics[width=0.7\linewidth]{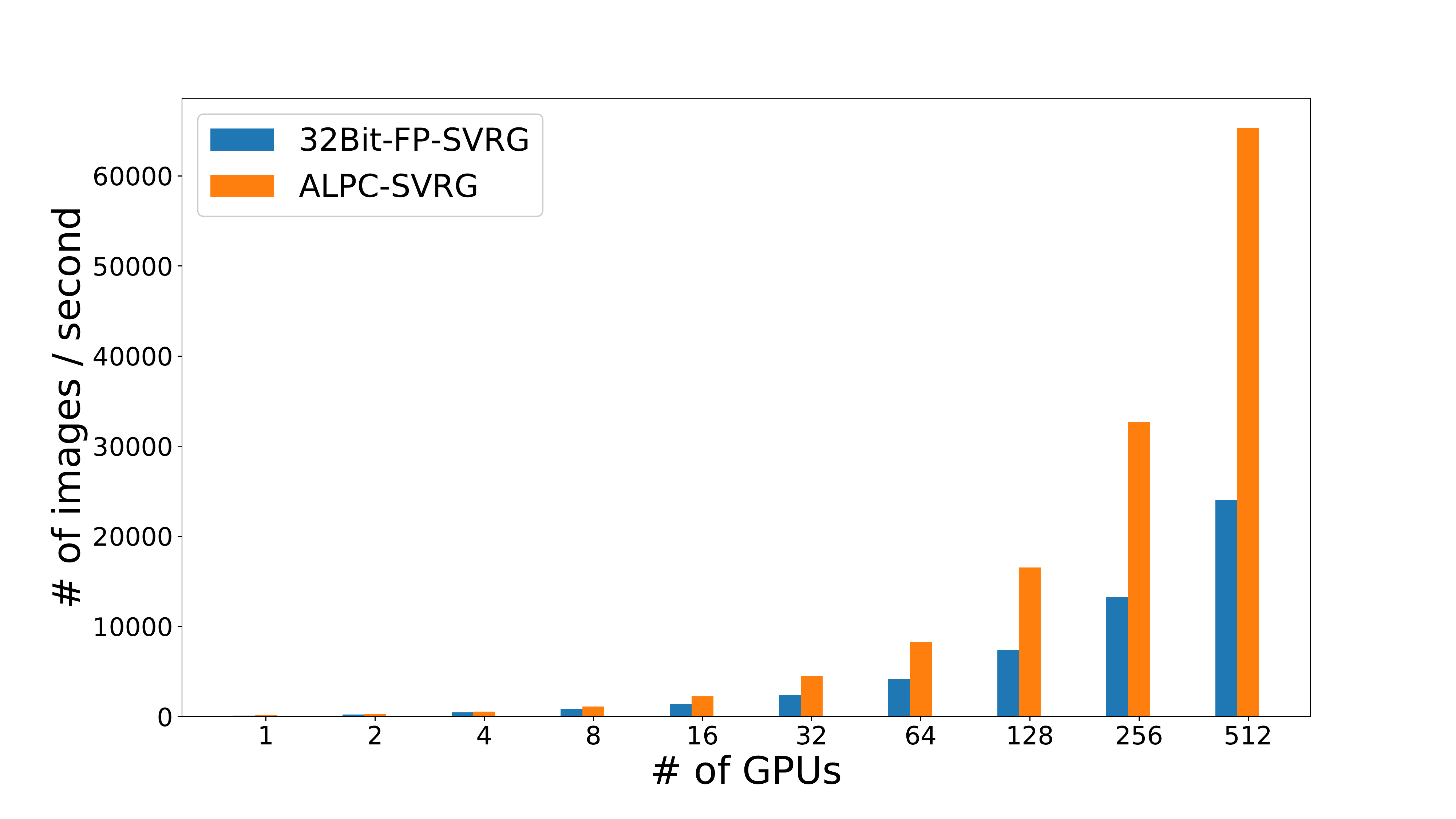}
\caption{Performance model: training throughput of ResNet-50 on ILSVRC-12 with the increasing of GPUs.}
\label{fig:preformance-model}
\end{figure}
\\
\noindent
\textbf{Performance Model.}
To verify the scalability of ALPC-SVRG,
we use the performance model proposed in \cite{yan2015performance}.
In this experiment, the performance denotes the training speed.
Similar to \cite{wu2018error,wen2017terngrad},  we use the lightweight profiling on the computation and communication time of a single machine to estimate the performance for larger distributed systems.
The major hardware statistics are as follows: Nvidia Tesla P40 GPU (8 units per node), Intel Xeon E5-2680 CPU and Mellanox ConnectX-3 Pro network card (with 40Gbps network connections).
\\
\\
\noindent
We setup experiments on ResNet-50 model on dataset ILSVRC-12 \cite{russakovsky2015imagenet}.
Figure~\ref{fig:preformance-model} presents the training throughput with the increasing of GPUs.
It shows that our algorithm can achieve significant speedup over full precision SVRG when applied to large-scale networks.
\vskip -0.1in
\section{CONCLUSION}
In this paper, we propose LPC-SVRG and its acceleration ALPC-SVRG to achieve both fast convergence and lower communication complexity.
We present a new quantization method with gradient clipping adopted and mathematically analyze its convergence.
%
With \emph{double sampling}, we are able to combine the gradients of both full-precision and low-precision and then achieve acceleration.
Our analysis covers general nonsmooth composite problems, and
shows the same asymptotic convergence rate can be attained with only $O(\log{\sqrt{d}})$ bits (compared to full-precision of $32$ bits).
The experiments on linear models and deep neural networks demonstrate the fast convergence and communication efficiency of our algorithms.

\subsubsection*{Acknowledgements}
The work by Yue Yu is supported in part by the National Natural Science Foundation of China Grants $61672316$.
We would like to thank Chaobing Song for the useful discussions and thank anonymous reviewers for their insightful suggestions.
\bibliography{sample_paper.bib}
\bibliographystyle{abbrv}

\clearpage

\section*{Supplementary Materials}
\section{Proof of Unbiased Quantization Variance}
\begin{lemma}
  \label{lemma:un-variance}
  If $x \in \mathbb{R}$ is in the convex hull of $\dom(\delta,b)$, then the quantization variance can be bounded as
  \begin{equation}
    \mathbf{E}[\big( Q_{(\delta,b)}(x) - x \big)^2] \leq \frac{\delta^2}{4}.
  \end{equation}
\end{lemma}
\noindent
\emph{Proof.}
From the manuscript we know that if $x$ is in the convex hull of $\dom(\delta,b)$, then it will be stochastically rounded up or down.
Without loss of generality, let $z+\delta$ and $z$ be the up and down quantization values respectively, then
\begin{equation*}
  Q_{(\delta,b)}(x) =
  \begin{cases}
    z \quad \quad \ \ \with-probability \ \frac{z+\delta-x}{\delta},\\
    z+\delta \quad   \otherwise.
  \end{cases}
\end{equation*}
Note that if $x$ equals to the smallest or largest value in $\dom(\delta,b)$,
%
then $Q_{(\delta,b)}(x) = x$ from the above definition of function $Q$.
Firstly, it can be verified that $\mathbf{E}\big[ Q_{(\delta,b)}(x) \big] = x$, then we have
\begin{equation}
  \begin{aligned}
  \mathbf{E}[\big( Q_{(\delta,b)}(x) - x \big)^2] = \frac{z+\delta-x}{\delta}(z-x)^2 + \frac{x-z}{\delta}(z+\delta-x)^2
=(x-z)(z+\delta-x) \leq \frac{\delta^2}{4}.
\end{aligned}
  \end{equation}
\QEDB

\section{Proof of Theorem $2$}
\begin{lemma} 
  \label{lemma:quantize-variance}
  For $\omega \in \mathbb{R}^d$, $\lambda \in (0,1]$, if $\delta = \frac{\lambda ||\omega||_{\infty}}{2^{b-1}-1}$,  then
  \begin{equation}
    \label{eq:quantization-variance}
  \mathbf{E}|| Q_{(\delta,b)}(\omega)- \omega||^2 \leq \frac{ (d-d_{\lambda}) \delta^2 }{4} + d_{\lambda}(1-\lambda)^2||\omega||^2,
\end{equation}
where $d_{\lambda}$ is the number of coordinates in $\omega$ exceeding $\dom(\delta, b)$.
\end{lemma}
\noindent
\emph{Proof.}
Since the squared norm $||Q_{(\delta,b)}(\omega) - \omega ||^2$ separates along dimensions, it suffices to consider a single coordinate $\omega_i$.
%
If $\omega_i$ is in the convex hull of $ \dom(\delta, b)$, then according to Lemma~\ref{lemma:un-variance} we have
\begin{equation}
  \mathbf{E}[\big( Q_{(\delta,b)}(\omega_i) - \omega_i \big)^2] \leq \frac{\delta^2}{4}.
\end{equation}
On the other hand, if $\omega_i$ is not in the convex hull of $\dom(\delta, b)$, then $Q_{(\delta,b)}(\omega_i)$ is either the smallest or the largest value of  $\dom(\delta,b)$.
Therefore,
 $\big( Q_{(\delta,b)}(\omega_i) - \omega_i \big)^2 \leq (\lambda||\omega||_{\infty} - ||\omega||_{\infty})^2 = (1-\lambda)^2||\omega||^2_{\infty} \leq (1-\lambda)^2||\omega||^2$ .
 Summing up over all dimensions we get (\ref{eq:quantization-variance}).
 \QEDB

\begin{lemma}
  \label{lemma:nonconvex-unbiased-variance} 
  For the iterates $x_t^{s+1}$, $\tilde{x}^s$ in Algorithm $1$, define $g_{t} \triangleq \frac{1}{B}\sum\limits_{j=1}^B \big[ \nabla f_j(x_{t}^{s+1}) - \nabla f_j(\tilde{x}^s) \big]+ \nabla f(\tilde{x}^s) $,
  where each element $j$ is uniformly and independently sampled from \{1,...,n\}, we have
  \begin{equation}
  \mathbf{E}||g_{t} - \nabla f(x_{t}^{s+1})||^2 \leq \frac{L^2}{B} \mathbf{E}||x_{t}^{s+1} - \tilde{x}^s||^2.
\end{equation}
\end{lemma}
\noindent \emph{Proof.}
\begin{equation}
\begin{aligned}
  \mathbf{E}||g_{t} - \nabla f(x_{t}^{s+1})||^2
  &= \mathbf{E} || \frac{1}{B}\sum\limits_{j=1}^B \big[ \nabla f_j(x_{t}^{s+1}) - \nabla f_j(\tilde{x}^s) + \nabla f(\tilde{x}^s) - \nabla f(x_{t}^{s+1})\big] ||^2\\
  &= \frac{1}{B^2} \sum\limits_{j=1}^B \mathbf{E} || \nabla f_j(x_{t}^{s+1}) - \nabla f_j(\tilde{x}^s) + \nabla f(\tilde{x}^s) - \nabla f(x_{t}^{s+1}) ||^2\\
  & \leq \frac{1}{B^2} \sum\limits_{j=1}^B \mathbf{E} || \nabla f_j(x_{t}^{s+1}) - \nabla f_j(\tilde{x}^s) ||^2\\
  & \leq \frac{L^2}{B} \mathbf{E} ||x_t^{s+1} - \tilde{x}^s||^2,
\end{aligned}
\end{equation}
where the first inequality uses $\mathbf{E}||x - \mathbf{E}x||^2 \leq \mathbf{E} ||x||^2$ and the last inequality follows from the Lipschitz smooth property of $f_j(x)$.
\QEDB

\begin{lemma}
  \label{lemma:nonconvex-biased-variance}
  Denote $d_{\lambda} = \max_i\{d_{\lambda}^i\}$, where
  $d_{\lambda}^i$ is the number of coordinates in $u_t^i$ exceeding $\dom(\delta_t^i,b)$.
Under communication scheme \textbf{(a)} or \textbf{(b)}, i.e., $\tilde{u}_t = \frac{1}{N}\sum\limits_{i=1}^N\tilde{u}_t^i$, we have
  \begin{equation}
    \label{eq:5}
    \mathbf{E}||v_t^{s+1} - \nabla f(x_t^{s+1})||^2 \leq 2L^2\Big[ \frac{d\lambda^2}{4(2^{b-1}-1)^2} + d_{\lambda}(1-\lambda)^2 + \frac{1}{NB}\Big]\mathbf{E}||x_t^{s+1} - \tilde{x}^s||^2.
  \end{equation}
  Note that $\delta_t^i$ has different value for \textbf{(a)} and \textbf{(b)}.
\end{lemma}
\noindent
\emph{Proof.}
\textbf{Case $1$.} First of all, we consider communication scheme \textbf{(a)}, i.e., $\tilde{u}_t = \frac{1}{N}\sum\limits_{i=1}^N\tilde{u}_t^i = \frac{1}{N}\sum\limits_{i=1}^NQ_{(\delta_t^i,b)}(u_t^i)$, where $\delta_t^i = \frac{\lambda ||u_t^i||_{\infty}}{2^{b-1}-1}$.
Therefore
\begin{equation}
  \begin{aligned}
    \label{eq:6}
  &\mathbf{E}||v_t^{s+1} - \nabla f(x_t^{s+1})||^2\\
  &= \mathbf{E} ||\frac{1}{N} \sum\limits_{i=1}^{N}Q_{(\delta_{t}^i, b)} (u_{t}^i) - \frac{1}{N}\sum\limits_{i=1}^N u_{t}^i + \frac{1}{N}\sum\limits_{i=1}^N u_{t}^i + \nabla f(\tilde{x}^s) - \nabla f(x_{t}^{s+1})||^2\\
  &\leq 2 \underbrace{\mathbf{E} ||\frac{1}{N} \sum\limits_{i=1}^{N}Q_{(\delta_{t}^i, b)} (u_{t}^i) - \frac{1}{N}\sum\limits_{i=1}^N u_{t}^i||^2 }_{L_1} + 2\mathbf{E}|| \frac{1}{N}\sum\limits_{i=1}^N u_{t}^i + \nabla f(\tilde{x}^s) - \nabla f(x_{t}^{s+1})||^2,
\end{aligned}
\end{equation}
where $L_1$ can be bounded as follows.
\begin{equation}
  \begin{aligned}
    \label{eq:7}
  L_1 &\leq \frac{1}{N}\sum\limits_{i=1}^N \mathbf{E} || Q_{(\delta_{t}^i, b)} (u_{t}^i) - u_{t}^i||^2\\
  &\leq \frac{1}{N}\sum\limits_{i=1}^N \mathbf{E} \Big[ \frac{d\lambda^2}{4(2^{b-1}-1)^2} + d_{\lambda}^i(1-\lambda)^2 \Big] ||u_t^i||^2\\
  &\leq \frac{1}{N}\Big[ \frac{d\lambda^2}{4(2^{b-1}-1)^2} + d_{\lambda}(1-\lambda)^2\Big] \sum\limits_{i=1}^N \mathbf{E}||u_{t}^i||^2\\
  &\leq L^2\Big[ \frac{d\lambda^2}{4(2^{b-1}-1)^2} + d_{\lambda}(1-\lambda)^2\Big]\mathbf{E}||x_t^{s+1} - \tilde{x}^s||^2,
\end{aligned}
\end{equation}
where the second inequality uses Lemma~\ref{lemma:quantize-variance} and the last inequality is due to the smoothness of $f_i(x)$.
Substituting (\ref{eq:7}) into (\ref{eq:6}) and using Lemma \ref{lemma:nonconvex-unbiased-variance}, we get (\ref{eq:5}).

\noindent
\textbf{Case $2$.}
If employing communication scheme \textbf{(b)}, we have $\tilde{u}_t = \frac{1}{N}\sum\limits_{i=1}^N\tilde{u}_t^i = \frac{1}{N}\sum\limits_{i=1}^N Q_{(\delta_t,b)}(u_t^i)$, where $\delta_t = \max_i\{\delta_t^i\}$.
Denote $j = \argmax_i{||u_t^i||_{\infty}}$, therefore $\delta_t = \frac{\lambda ||u_t^j||_{\infty}}{2^{b-1}-1}$.
Then let $\delta_t^i = \delta_t$ for all $i$, putting back into (\ref{eq:6}) we obtain
\begin{equation}
  \begin{aligned}
    &\mathbf{E}||v_t^{s+1} - \nabla f(x_t^{s+1})||^2\\
    &\leq 2 \underbrace{\mathbf{E} ||\frac{1}{N} \sum\limits_{i=1}^{N}Q_{(\delta_{t}, b)} (u_{t}^i) - \frac{1}{N}\sum\limits_{i=1}^N u_{t}^i||^2 }_{L_1^{\prime}} + 2\mathbf{E}|| \frac{1}{N}\sum\limits_{i=1}^N u_{t}^i + \nabla f(\tilde{x}^s) - \nabla f(x_{t}^{s+1})||^2,
  \end{aligned}
\end{equation}
where $L_1^{\prime}$ has a bound of
\begin{equation}
  \begin{aligned}
    \label{eq:prime}
  L_1^{\prime} &\leq \frac{1}{N}\sum\limits_{i=1}^N \mathbf{E} || Q_{(\delta_{t}, b)} (u_{t}^i) - u_{t}^i||^2\\
  &\leq \frac{1}{N}\sum\limits_{i=1}^N \mathbf{E} \Big[ \frac{d\lambda^2}{4(2^{b-1}-1)^2} + d_{\lambda}^i(1-\lambda)^2 \Big] ||u_t^j||^2\\
  &\leq \frac{1}{N}\Big[ \frac{d\lambda^2}{4(2^{b-1}-1)^2} + d_{\lambda}(1-\lambda)^2\Big] \sum\limits_{i=1}^N \mathbf{E}||u_{t}^j||^2\\
  &\leq L^2\Big[ \frac{d\lambda^2}{4(2^{b-1}-1)^2} + d_{\lambda}(1-\lambda)^2\Big]\mathbf{E}||x_t^{s+1} - \tilde{x}^s||^2.
\end{aligned}
\end{equation}
We adopt
\begin{equation}
  \mathbf{E} || Q_{(\delta_{t}, b)} (u_{t}^i) - u_{t}^i||^2 \leq \mathbf{E} \Big[ \frac{d\lambda^2}{4(2^{b-1}-1)^2} + d_{\lambda}^i(1-\lambda)^2 \Big] ||u_t^j||^2
\end{equation}
 in the second inequality in (\ref{eq:prime}), which can be verified following the proof of Lemma~\ref{lemma:quantize-variance} with $\delta_t = \frac{\lambda ||u_t^j||_{\infty}}{2^{b-1}-1}$.
Putting the above inequalities together we obtain (\ref{eq:5}).
\QEDB

\begin{lemma}
  \label{lemma:nonconvex-aggra-biased-variance}
  Under communication scheme \textbf{(c)} with $\tilde{u}_t = Q_{(\delta_{t}, b)}(\frac{1}{N}\sum\limits_{i=1}^N \tilde{u}_t^i)$, $\tilde{u}_t^i = Q_{(\delta_t,b)}(u_t^i)$,
$\delta_t = \max_i\{\delta_t^i\}$,
  we obtain
  \begin{equation}
    \label{eq:nonconvex-aggra-biased-variance}
  \mathbf{E}||v_t^{s+1} - \nabla f(x_t^{s+1})||^2 \leq 2L^2\Big[ \frac{3d\lambda^2}{8(2^{b-1}-1)^2} + d_{\lambda}(1-\lambda)^2 + \frac{1}{NB} \Big]\mathbf{E}||x_t^{s+1} - \tilde{x}^s||^2,
\end{equation}
where $d_{\lambda} = \max_i\{d_{\lambda}^i\}$,
$d_{\lambda}^i$ is the number of coordinates in $u_t^i$ exceeding $\dom(\delta_t^i,b)$.
\end{lemma}
\noindent
\emph{Proof.}
\begin{equation}
  \begin{aligned}
\mathbf{E}||v_t^{s+1} - \nabla f(x_t^{s+1})||^2 =  \underbrace{\mathbf{E}||Q_{(\delta_t,b)}(\frac{1}{N}\sum\limits_{i=1}^N \tilde{u}_t^i) - \frac{1}{N}\sum\limits_{i=1}^N \tilde{u}_t^i ||^2}_{L_2} + \underbrace{\mathbf{E}||\frac{1}{N}\sum\limits_{i=1}^N \tilde{u}_t^i + \nabla f(\tilde{x}^s) - \nabla f(x_t^{s+1})||^2}_{L_3}\\
\end{aligned}
\end{equation}
where the equality holds
because $Q_{(\delta_t,b)}(\frac{1}{N}\sum\limits_{i=1}^N \tilde{u}_t^i)$ is an unbiased quantization
(note that each $u_t^i$ is quantized using $\delta_t$, therefore, all coordinates of $\frac{1}{N}\sum\limits_{i=1}^N \tilde{u}_t^i$ are in the convex hull of  $\dom(\delta_t,b)$).

\noindent
From  Lemma~\ref{lemma:un-variance} and \textbf{Case 2} in
Lemma \ref{lemma:nonconvex-biased-variance}
we obtain
\begin{equation}
  \begin{aligned}
L_2 & \leq \frac{d\delta_t^2}{4} \\
& \leq \frac{d\lambda^2||u_t^j||^2_{\infty}}{4(2^{b-1}-1)^2} ,\quad j=\argmax_i{||u_t^i||_{\infty}}\\
& \leq \frac{L^2d\lambda^2}{4(2^{b-1}-1)^2}\mathbf{E}||x_t^{s+1}-\tilde{x}^s||^2
\end{aligned}
\end{equation}
and
\begin{equation}
  L_3 \leq 2L^2\Big[ \frac{d\lambda^2}{4(2^{b-1}-1)^2} + d_{\lambda}(1-\lambda)^2 + \frac{1}{NB}\Big]\mathbf{E}||x_t^{s+1} - \tilde{x}^s||^2.
  \end{equation}
  Putting them together, we get (\ref{eq:nonconvex-aggra-biased-variance}).
  \QEDB
\\
\\
\noindent
\textbf{\emph{Proof of Theorem $2$.}}
Define $\bar{x}_{t+1}^{s+1} = \prox_{\eta h}(x_t^{s+1} - \eta \nabla f(x_t^{s+1}))$.
Following the proof of Theorem $5$ in \cite{Reddi2016Fast} (equations $(8)$-$(12)$), we get
\begin{equation}
  \mathbf{E} \Big[P(x_{t+1}^{s+1})\Big] \leq \mathbf{E} \Big[ P(x_t^{s+1}) + \frac{\eta}{2}||v_t^{s+1} - \nabla f(x_t^{s+1})||^2 + (L-\frac{1}{2\eta})||\bar{x}_{t+1}^{s+1} - x_t^{s+1}||^2 + (\frac{L}{2} - \frac{1}{2\eta})||x_{t+1}^{s+1} - x_t^{s+1}||^2\Big].
\end{equation}
If adopting communication scheme \textbf{(a)} or \textbf{(b)}, combining Lemma \ref{lemma:nonconvex-biased-variance}, we have
\begin{equation}
  \begin{aligned}
  \mathbf{E} \Big[P(x_{t+1}^{s+1}) \Big]
  \leq & \mathbf{E} \Big[ P(x_t^{s+1}) + \eta L^2\Big[  \frac{d\lambda^2}{4(2^{b-1}-1)^2} + d_{\lambda}(1-\lambda)^2 + \frac{1}{NB} \Big]||x_t^{s+1} - \tilde{x}^s||^2 \\
  &+ (L-\frac{1}{2\eta})||\bar{x}_{t+1}^{s+1} - x_t^{s+1}||^2 + (\frac{L}{2} - \frac{1}{2\eta})||x_{t+1}^{s+1} - x_t^{s+1}||^2\Big].
\end{aligned}
\end{equation}
Define $R_t^{s+1} \triangleq \mathbf{E}\Big[P(x_t^{s+1}) + c_t||x_t^{s+1} - \tilde{x}^s||^2\Big]$
 and a sequence $\{c_t\}_{t=0}^{m}$ with $c_m = 0$ and
 $c_t = c_{t+1}(1+\beta) + \eta L^2 \Big[ \frac{d\lambda^2}{4(2^{b-1}-1)^2} + d_{\lambda}(1-\lambda)^2 + \frac{1}{NB} \Big]$, where $\beta = \frac{1}{m}$.
 Therefore $\{c_t\}$ is a decreasing sequence. We first derive the bound of $c_0$ in the following.
 \begin{equation}
   \begin{aligned}
     \label{eq:bound-c}
     c_0 &\leq \eta L^2 \Big[ \frac{d\lambda^2}{4(2^{b-1}-1)^2} + d_{\lambda}(1-\lambda)^2 + \frac{1}{NB} \Big] \cdot \frac{(1+\beta)^m - 1}{\beta}\\
     &\leq 2m\eta L^2 \Big[ \frac{d\lambda^2}{4(2^{b-1}-1)^2} + d_{\lambda}(1-\lambda)^2 + \frac{1}{NB} \Big]
   \end{aligned}
 \end{equation}
 where the second inequality uses $\beta = \frac{1}{m}$. Denote $\eta = \frac{\rho}{L}$.
 Then inequality (\ref{eq:bound-c}) can be simplified as
 \begin{equation}
   \label{eq:c0}
   c_0 \leq 2m\rho L\Big[  \frac{d\lambda^2}{4(2^{b-1}-1)^2} + d_{\lambda}(1-\lambda)^2 + \frac{1}{NB}\Big].
 \end{equation}
On the other hand,
\begin{equation}
  \begin{aligned}
    R_{t+1}^{s+1} &= \mathbf{E} \Big[ P(x_{t+1}^{s+1}) + c_{t+1}||x_{t+1}^{s+1} - \tilde{x}^s||^2 \Big]\\
    &\leq \mathbf{E} \Big[ P(x_{t+1}^{s+1}) + c_{t+1}(1+\frac{1}{\beta})||x_{t+1}^{s+1} - x_t^{s+1}||^2 + c_{t+1}(1+\beta)||x_t^{s+1}-\tilde{x}^s||^2\Big]\\
    &\leq \mathbf{E}\Big[ P(x_t^{s+1}) + c_t||x_t^{s+1} - \tilde{x}^s||^2 + (c_{t+1}(1+\frac{1}{\beta}) + \frac{L}{2} - \frac{1}{2\eta})||x_{t+1}^{s+1} - x_t^{s+1}||^2 + (L-\frac{1}{2\eta})||\bar{x}_{t+1}^{s+1} - x_t^{s+1}||^2\Big].
  \end{aligned}
\end{equation}
Now we derive the bound for $\rho$ and $b$ to make sure $(c_{t+1}(1+\frac{1}{\beta}) + \frac{L}{2} - \frac{1}{2\eta}) \leq 0$,
and it suffices to let $c_0(1+\frac{1}{\beta}) + \frac{L}{2} \leq \frac{1}{2\eta}$.
Combining (\ref{eq:c0}) and $\beta = \frac{1}{m}$, $\eta = \frac{\rho}{L}$, we only need to guarantee
\begin{equation}
  8m^2\rho^2\Big[ \frac{d\lambda^2}{4(2^{b-1}-1)^2} + d_{\lambda}(1-\lambda)^2 + \frac{1}{NB} \Big] + \rho \leq 1.
\end{equation}
If the above constraint holds, then
\begin{equation}
  \begin{aligned}
    R_{t+1}^{s+1} &\leq R_t^{s+1} + (L-\frac{1}{2\eta})\mathbf{E}||\bar{x}_{t+1}^{s+1} - x_t^{s+1}||^2\Big].
  \end{aligned}
\end{equation}
Summing it up over $t=0$ to $m-1$ and $s=0$ to $S-1$, using $c_m = 0$, $x_0^{s+1} = \tilde{x}^s$ and $x_m^{s+1} = \tilde{x}^{s+1}$ we get
\begin{equation}
  (\frac{1}{2\eta} - L) \sum\limits_{s=0}^{S-1}\sum\limits_{t=0}^{m-1}\mathbf{E}||\bar{x}_{t+1}^{s+1} - x_t^{s+1}||^2 \leq P(x^0) - P(x^*).
\end{equation}
Applying the definition of $G_\eta(x_t^{s+1})$, we obtain results in Theorem $2$.
Moreover, the analysis of communication scheme \textbf{(c)} can be similarly obtained using the above proof steps.
\QEDB

\section{Proof of ALPC-SVRG}
\begin{lemma}
  \label{lemma:unbiased-variance}
  For $\hat{v}_{k+1} = \frac{1}{B}\sum\limits_{j=1}^B \big[ \nabla f_j(x_{k+1}) - \nabla f_j(\tilde{x}^s) \big]+ \nabla f(\tilde{x}^s) $, where each element $j$ is uniformly and independently sampled from \{1,...,n\}, we have
  \begin{equation}
  \mathbf{E}||\hat{v}_{k+1} - \nabla f(x_{k+1})||^2 \leq \frac{2L}{B} \mathbf{E}\big[ f(\tilde{x}^s) - f(x_{k+1}) - \langle \nabla f(x_{k+1}), \tilde{x}^{s} - x_{k+1}\rangle\big].
\end{equation}
\end{lemma}
\noindent \emph{Proof.}
\begin{equation}
\begin{aligned}
  \mathbf{E}||\hat{v}_{k+1} - \nabla f(x_{k+1})||^2
  & \leq \frac{1}{B^2} \sum\limits_{j=1}^B \mathbf{E} || \nabla f_j(x_{k+1}) - \nabla f_j(\tilde{x}^s) ||^2\\
  & \leq \frac{2L}{B} \mathbf{E} \big[ f(\tilde{x}^s) - f(x_{k+1}) - \langle \nabla f(x_{k+1}), \tilde{x}^s - x_{k+1} \rangle\big],
\end{aligned}
\end{equation}
where the first inequality follows from Lemma \ref{lemma:nonconvex-unbiased-variance} and the last inequality adopts the Lipschitz smooth property of $f_j(x)$.
\QEDB

\begin{lemma}
  \label{lemma:2-biased-variance}
  Denote $d_{\lambda} = \max_i\{d_{\lambda}^i\}$, where
  $d_{\lambda}^i$ is the number of coordinates in $u_{k+1}^i$ exceeding $\dom(\delta_{k+1}^i,b)$, then we have
  \begin{equation}
    \label{eq:biased-variance}
    \mathbf{E}||v_{k+1} - \nabla f(x_{k+1})||^2 \leq 4L\Big[ \frac{d\lambda^2}{4(2^{b-1}-1)^2} + d_{\lambda}(1-\lambda)^2 + \frac{1}{NB}\Big] \mathbf{E}\big[ f(\tilde{x}^s) - f(x_{k+1}) - \langle \nabla f(x_{k+1}), \tilde{x}^s - x_{k+1}\rangle \big].
  \end{equation}
\end{lemma}
\noindent
\emph{Proof.}
\begin{equation}
  \begin{aligned}
&\mathbf{E} ||v_{k+1} - \nabla f(x_{k+1})||^2 \\
&\leq 2 \underbrace{\mathbf{E} ||\frac{1}{N} \sum\limits_{i=1}^{N}Q_{(\delta_{k+1}^i, b)} (u_{k+1}^i) - \frac{1}{N}\sum\limits_{i=1}^N u_{k+1}^i||^2 }_{A_1}+ 2 \underbrace{\mathbf{E}|| \frac{1}{N}\sum\limits_{i=1}^N u_{k+1}^i + \nabla f(\tilde{x}^s) - \nabla f(x_{k+1})||^2}_{A_2}.
\end{aligned}
\end{equation}
Using the same arguments of Lemma~\ref{lemma:unbiased-variance} we obtain
\begin{equation}
  A_2 \leq \frac{2L}{NB}\mathbf{E}\big[ f(\tilde{x}^s) - f(x_{k+1}) - \langle \nabla f(x_{k+1}), \tilde{x}^s - x_{k+1}\rangle \big].
\end{equation}
Moreover,
\begin{equation}
  \begin{aligned}
    A_1 &\leq \frac{1}{N}\sum\limits_{i=1}^N \mathbf{E} || Q_{(\delta_{k+1}^i, b)} (u_{k+1}^i) - u_{k+1}^i||^2\\
    &\leq \frac{1}{N}\sum\limits_{i=1}^N \mathbf{E} \Big[ \frac{d\lambda^2}{4(2^{b-1}-1)^2} + d_{\lambda}^i(1-\lambda)^2\Big] ||u_{k+1}^i||^2\\
    &\leq \frac{1}{N}\Big[ \frac{d\lambda^2}{4(2^{b-1}-1)^2} + d_{\lambda}(1-\lambda)^2\Big] \sum\limits_{i=1}^N \mathbf{E}||u_{k+1}^i||^2\\
    &\leq 2L\Big[ \frac{d\lambda^2}{4(2^{b-1}-1)^2} + d_{\lambda}(1-\lambda)^2\Big]\mathbf{E}\big[ f(\tilde{x}^s) - f(x_{k+1}) - \langle \nabla f(x_{k+1}), \tilde{x}^s - x_{k+1}\rangle \big],
  \end{aligned}
\end{equation}
where the second inequality follows from Lemma \ref{lemma:quantize-variance}.
Putting them together, we obtain (\ref{eq:biased-variance}).
\QEDB

\begin{lemma}
  \label{lemma:biased-unbiased-variance}
  \begin{equation}
    \label{eq:biased-unbiased-variance}
  \mathbf{E}||\hat{v}_{k+1} - v_{k+1}||^2 \leq 4L\Big[ \frac{1}{2B} + \frac{d\lambda^2}{4(2^{b-1}-1)^2} + d_{\lambda}(1-\lambda)^2 + \frac{1}{NB}\Big] \mathbf{E}\big[ f(\tilde{x}^s) - f(x_{k+1}) - \langle \nabla f(x_{k+1}), \tilde{x}^s - x_{k+1}\rangle \big].
\end{equation}
\end{lemma}
\noindent
\emph{Proof.}
\begin{equation}
\begin{aligned}
\mathbf{E}||\hat{v}_{k+1} - v_{k+1}||^2
&= \mathbf{E} ||\hat{v}_{k+1} - \nabla f(x_{k+1}) + \nabla f(x_{k+1}) - v_{k+1}||^2\\
&= \mathbf{E}||\hat{v}_{k+1} - \nabla f(x_{k+1}) ||^2 + {\mathbf{E} ||v_{k+1} - \nabla f(x_{k+1}) ||^2}.
\end{aligned}
\end{equation}
The second equality holds because the mini-batches for calculating $\hat{v}_{k+1}$ and $v_{k+1}$ are independent and $\mathbf{E}\hat{v}_{k+1} = \nabla f(x_{k+1})$.
Combining Lemma~\ref{lemma:unbiased-variance} and Lemma~\ref{lemma:2-biased-variance}, we obtain (\ref{eq:biased-unbiased-variance}).
\QEDB

\begin{lemma}
  \label{lemma:update-y}
  Define
  \begin{equation}
\Prog (x_{k+1}) \triangleq - min_{y} \{2L||y-x_{k+1}||^2 + \langle v_{k+1}, y-x_{k+1}\rangle + h(y) - h(x_{k+1}) \},
  \end{equation}
  then from the update rule of $y$, we obtain
  \begin{equation}
    \label{eq:y-inequa}
    \mathbf{E} \big[ P(x_{k+1}) - P(y_{k+1})\big] \geq \mathbf{E} \big[\Prog (x_{k+1}) - \frac{1}{6L}||\nabla f(x_{k+1}) - v_{k+1}||^2 \big].
   \end{equation}
\end{lemma}
\noindent
\emph{Proof Sketch.} (\ref{eq:y-inequa}) follows from the proof of Lemma $3.3$ in \cite{allen2017katyusha} with different coefficients.
\QEDB

\begin{lemma} If $h(x)$ is $\sigma$-strongly convex, then for any $u \in \mathbb{R}^d$, we have (Lemma $3.5$ in \cite{allen2017katyusha})
\label{lemma:z-inequa}
  \begin{equation}
  \alpha \langle \hat{v}_{k+1}, z_{k+1} - u \rangle + \alpha h(z_{k+1}) - \alpha h(u) \leq -\frac{1}{2}||z_k - z_{k+1}||^2 + \frac{1}{2}||z_k - u||^2 - \frac{1+\alpha \sigma}{2}||z_{k+1}-u||^2.
\end{equation}
\end{lemma}
\noindent

\begin{lemma}
  \label{lemma:middle-1}
  Let $\tau_2 = \frac{5}{3}\zeta + \frac{1}{2B}$, $\zeta = \frac{d\lambda^2}{4(2^{b-1}-1)^2} + d_{\lambda}(1-\lambda)^2 + \frac{1}{NB}$, $\alpha = \frac{1}{6\tau_1 L}$.
  Suppose with proper choice of parameters $B,b,\lambda$, we have $\tau_2 \leq \frac{1}{2}$, then
  \begin{equation}
    \label{eq:middle-1}
    \begin{aligned}
  & \mathbf{E} \Big[ \alpha \langle \nabla f(x_{k+1}), z_k - u\rangle - \alpha h(u) \Big]\\
  &\leq \mathbf{E} \Big[ \frac{\alpha}{\tau_1}\big[ P(x_{k+1}) - P(y_{k+1}) + \tau_2 P(\tilde{x}^s) - \tau_2 f(x_{k+1}) - \tau_2 \langle \nabla f(x_{k+1}), \tilde{x}^s - x_{k+1} \rangle\big] \\
  &+ \frac{\alpha}{\tau_1} (1-\tau_1-\tau_2) h(y_k) - \frac{\alpha}{\tau_1}h(x_{k+1})
  + \frac{1}{2}||z_k - u||^2 - \frac{1+\alpha \sigma}{2}||z_{k+1}-u||^2 \Big].
    \end{aligned}
  \end{equation}
\end{lemma}
\noindent
\emph{Proof.}
\begin{equation}
  \begin{aligned}
    \label{eq:13}
  & \mathbf{E} \Big[ \alpha \langle \hat{v}_{k+1}, z_k - u \rangle + \alpha h(z_{k+1}) - \alpha h(u) \Big] \\
  = & \mathbf{E} \Big[ \alpha \langle v_{k+1}, z_k - z_{k+1}\rangle + \alpha \langle \hat{v}_{k+1} - v_{k+1}, z_k - z_{k+1} \rangle + \alpha \langle \hat{v}_{k+1},z_{k+1}-u \rangle + \alpha h(z_{k+1}) - \alpha h(u) \Big] \\
  \leq & \mathbf{E} \Big[ \alpha \langle v_{k+1}, z_k - z_{k+1}\rangle + \frac{\alpha}{\tau_1} \frac{1}{4L}||\hat{v}_{k+1} - v_{k+1}||^2 + \frac{1}{6}||z_k - z_{k+1}||^2 + \alpha \langle \hat{v}_{k+1},z_{k+1}-u \rangle + \alpha h(z_{k+1}) - \alpha h(u) \Big] \\
  \leq & \mathbf{E} \Big[ \alpha \langle v_{k+1}, z_k - z_{k+1}\rangle + \frac{\alpha}{\tau_1} \frac{1}{4L}||\hat{v}_{k+1} - v_{k+1}||^2 - \frac{1}{3}||z_k - z_{k+1}||^2 + \frac{1}{2}||z_k - u||^2 - \frac{1+\alpha \sigma}{2}||z_{k+1}-u||^2 \Big] \\
  \leq & \mathbf{E} \Big[ \alpha \langle v_{k+1}, z_k - z_{k+1}\rangle + \frac{\alpha}{\tau_1} (\zeta+\frac{1}{2B}) \mathbf{E}\big[ f(\tilde{x}^s) - f(x_{k+1}) - \langle \nabla f(x_{k+1}), \tilde{x}^{s} - x_{k+1}\rangle\big] - \frac{1}{3}||z_k - z_{k+1}||^2  \\
  &+  \frac{1}{2}||z_k - u||^2 - \frac{1+\alpha \sigma}{2}||z_{k+1}-u||^2\Big],
\end{aligned}
\end{equation}
where the first inequality uses Young's inequality and $\alpha = \frac{1}{6\tau_1 L}$,
the last two inequalities follow from Lemma~\ref{lemma:z-inequa} and Lemma~\ref{lemma:biased-unbiased-variance} respectively.
\noindent
Define $v \triangleq \tau_1 z_{k+1} + \tau_2 \tilde{x}^{s} + (1-\tau_1-\tau_2)y_k$, therefore $x_{k+1} - v = \tau_1(z_k-z_{k+1})$, then we obtain
\begin{equation}
  \begin{aligned}
    \label{eq:14}
  & \mathbf{E} \Big[\alpha \langle v_{k+1}, z_k - z_{k+1} \rangle - \frac{1}{3}||z_k - z_{k+1}||^2 \Big]\\
  = &\mathbf{E} \Big[ \frac{\alpha}{\tau_1}\langle v_{k+1}, x_{k+1} - v \rangle - \frac{1}{3\tau_1^2}||x_{k+1}-v||^2\Big]\\
  = & \mathbf{E} \Big[ \frac{\alpha}{\tau_1} \Big( \langle v_{k+1}, x_{k+1} - v \rangle - \frac{1}{3\alpha\tau_1}||x_{k+1}-v||^2 - h(v) + h(x_{k+1})\Big) + \frac{\alpha}{\tau_1}\Big(h(v) - h(x_{k+1})\Big)\Big]\\
  = & \mathbf{E} \Big[ \frac{\alpha}{\tau_1} \Big( \langle v_{k+1}, x_{k+1} - v \rangle - 2L||x_{k+1}-v||^2 - h(v) + h(x_{k+1})\Big) + \frac{\alpha}{\tau_1}\Big(h(v) - h(x_{k+1})\Big)\Big]\\
  \leq & \mathbf{E}\Big[ \frac{\alpha}{\tau_1}\Big(P(x_{k+1}) - P(y_{k+1}) + \frac{1}{6L}||v_{k+1} - \nabla f(x_{k+1})||^2\Big) + \frac{\alpha}{\tau_1}\Big(h(v) - h(x_{k+1})\Big)\Big]\\
\leq &\frac{\alpha}{\tau_1} \mathbf{E}\Big[ P(x_{k+1}) - P(y_{k+1})\Big] + \frac{\alpha}{\tau_1}\frac{2}{3}\zeta \mathbf{E}\Big[ f(\tilde{x}^s) - f(x_{k+1}) - \langle \nabla f(x_{k+1}), \tilde{x}^{s} - x_{k+1}\rangle\Big] + \frac{\alpha}{\tau_1}\mathbf{E} \Big[ h(v) - h(x_{k+1})\Big],
\end{aligned}
\end{equation}
where the third equality uses $\alpha = \frac{1}{6\tau_1 L}$, the first inequality follows from Lemma~\ref{lemma:update-y} and the last inequality adopts Lemma~\ref{lemma:2-biased-variance}.
Substituting (\ref{eq:14}) into (\ref{eq:13}) we get
\begin{equation}
  \begin{aligned}
& \mathbf{E} \Big[\alpha \langle  \hat{v}_{k+1}, z_k - u\rangle + \alpha h(z_{k+1})- \alpha h(u)\Big] \\
&\leq \mathbf{E} \Big[ \frac{\alpha}{\tau_1}\big[ P(x_{k+1}) - P(y_{k+1}) \big]+ \frac{\alpha}{\tau_1} \tau_2\Big[ f(\tilde{x}^s) - f(x_{k+1}) - \langle \nabla f(x_{k+1}), \tilde{x}^s - x_{k+1} \rangle\Big]\\
&+ \frac{1}{2}||z_k - u||^2 - \frac{1+\alpha \sigma}{2}||z_{k+1}-u||^2
+ \frac{\alpha}{\tau_1}\big[ \tau_1 h(z_{k+1}) + \tau_2 h(\tilde{x}^s) + (1-\tau_1-\tau_2)h(y_k) - h(x_{k+1})\big] \Big].
\end{aligned}
\end{equation}
Because $\hat{v}_{k+1}$ is unbiased, we get (\ref{eq:middle-1}) after rearranging terms.
\QEDB
\\
\\
\noindent
\textbf{\emph{Proof of Theorem $3$.}}
Starting form Lemma~\ref{lemma:middle-1}, following the proof of (\cite{allen2017katyusha}, Lemma $3.7$, Theorem $3.1$), we get
\begin{equation}
  \begin{aligned}
  &  \mathbf{E} \Big[ \frac{\tau_1+\tau_2-(1-\frac{1}{\theta})}{\tau_1}\theta \tilde{D}^{s+1} \cdot \sum\limits_{t=0}^{m-1}\theta^t \Big] \\
  &  \leq \mathbf{E} \Big[ \frac{1-\tau_1-\tau_2}{\tau_1}(D_{sm} - \theta^m D_{(s+1)m}) + \frac{\tau_2}{\tau_1} \tilde{D}^s \sum\limits_{t=0}^{m-1}\theta^t
    + \frac{1}{2\alpha}||z_{sm} - x^*||^2 - \frac{\theta^m}{2\alpha}||z_{(s+1)m} - x^*||^2 \Big],
  \end{aligned}
\end{equation}
where $\theta = (1+\alpha\sigma)$, $D_k \triangleq P(y_k) - P(x^*), \tilde{D}^s \triangleq P(\tilde{x}^s) - P(x^*)$.
\\
\\
\noindent
If $\frac{m\sigma}{L} \leq \frac{3}{2}$, then $\sqrt{\frac{m\sigma}{6L}} \leq \frac{1}{2}$.
Choosing $\alpha = \frac{1}{\sqrt{6m\sigma L}}$, then $\tau_1 = \frac{1}{6\alpha L} = m\sigma\alpha = \sqrt{\frac{m\sigma}{6L}} \leq \frac{1}{2}$ and $\alpha\sigma \leq \frac{1}{2m}$.
It can be verified that the above parameter settings guarantee \textbf{Case $1.$} in ( \cite{allen2017katyusha}, Theorem $1$),
 therefore, with the same arguments we arrive at
\begin{equation}
  \begin{aligned}
    \mathbf{E} \big[ P(\tilde{x}^S) - P(x^*) \big] \leq O((1+\alpha\sigma)^{-Sm}) \big[ P(x_0) - P(x^*) \big].
\end{aligned}
\end{equation}
\noindent
\QEDB
\\
\noindent
\textbf{\emph{Proof Sketch of Theorem $4$.}}
Let $\alpha_s = \frac{1}{6L\tau_{1,s}}, \tau_{1,s} = \frac{2}{s+4}$ and $\tau_2$ unchanged.
It can be verified that Lemma~\ref{lemma:middle-1} also holds in the current parameter setting (with $\sigma = 0$),
then plug Lemma~\ref{lemma:middle-1} into the proof of Theorem $4.1$ in \cite{allen2017katyusha}, we get
\begin{equation}
\begin{aligned}
\mathbf{E} \Big[ P(\tilde{x}^S) - P(x^*)\Big] \leq O(\frac{1}{mS^2}) \Big[ m(P(x_0) - P(x^*)) + L ||x_0 - x^*||^2 \Big].
\end{aligned}
\end{equation}
\QEDB

\end{document}